\documentclass[reqno]{amsart}
\usepackage{amsmath}
\allowdisplaybreaks[3]
\numberwithin{equation}{section}
\numberwithin{equation}{section}

\def\proof{\indent{\em Proof.\quad}}
\def\endproof{\hfill\hbox{$\sqcup$}\llap{\hbox{$\sqcap$}}\medskip}
\newtheorem{thm}{{\indent\bf Theorem}}[section]
\newtheorem{prop}{{\indent\bf Proposition}}[section]
\newtheorem{dfn}{{\indent\bf Definition}}[section]
\newtheorem{rmk}{{\indent\bf Remark}}[section]
\newtheorem{lem}{{\indent\bf Lemma}}[section]
\newtheorem{cor}{{\indent\bf Corollary}}[section]
\newtheorem{expl}{{\indent\bf Example}}[section]

\newcommand{\mb}{\mbox}
\newcommand{\hs}{\hspace}

\newcommand{\ol}{\overline}

\newcommand{\strl}{\stackrel}

\newcommand{\td}{\tilde}

\newcommand{\fr}{\frac}
\newcommand{\ed}{{\rm End}}
\newcommand{\edd}{\end{document}}
\newcommand{\be}{\begin{equation}}
\newcommand{\ee}{\end{equation}}

\newcommand{\lmx}{\left(\begin{matrix}}
\newcommand{\rmx}{\end{matrix}\right)}
\newcommand{\ldt}{\left|\begin{matrix}}
\newcommand{\rdt}{\end{matrix}\right|}

\newcommand{\tr}{{\rm tr\,}}

\newcommand{\veps}{\varepsilon}
\newcommand{\bbr}{{\mathbb R}}

\newcommand{\bbc}{{\mathbb C}}

\newcommand{\ba}{\begin{array}}
\newcommand{\ea}{\end{array}}
\newcommand{\nnm}{\nonumber}
\newcommand{\beal}{\begin{align}}
\newcommand{\eal}{\end{align}}
\newcommand{\bea}{\begin{eqnarray}}
\newcommand{\eea}{\end{eqnarray}}

\newcommand{\pp}[2]{\fr{\partial #1}{\partial #2}}
\newcommand{\ppp}[3]{\fr{\partial^2 #1}{\partial #2\partial #3}}
\newcommand{\dd}[2]{\fr{d #1}{d #2}}
\begin{document}

\title[On the Calabi composition of multiple affine hyperspheres]
{On the Calabi composition of multiple affine hyperspheres}%
\author{Xingxiao Li}%

%\date{}%
%\dedicatory{}%
%\commby{}%
% ----------------------------------------------------------------
\begin{abstract}
In this paper, we explicitly construct the Calabi composition of multiple affine hyperspheres possibly including some points viewing as $0$-dimensional hypersheres. Then we compute all the basic affine invariants of the composed affine hyperspheres, proving that the composed affine hypersphere is symmetric one if and only if each of its composing factors of positive dimension is symmetric.
\end{abstract}
\maketitle
% ----------------------------------------------------------------

\tableofcontents

\section{Introduction}

As we know, affine hyperspheres are the most important objects in affine differential geometry of nondegenerate hypersursurfaces, drawing great attention of many geometers.
In fact, affine hyperspheres seems simple in definition but they do form a very large class of hypersurfaces, the study of which is fruitful in recent twenty years. See for example, the proof of the Calabi's conjecture (\cite{amli90}, \cite{amli92}), the classification of hyperspheres of constant affine curvatures (\cite{vra-li-sim91},  \cite{wang93}, \cite{kri-vra99}), and in \cite{hu-li-vra11} the complete classification of locally strongly convex hypersurfaces with parallel Fubini-Pick forms which form a special class of hyperbolic affine hyperspheres (for some special cases, see \cite{dil-vra-yap94}, \cite{hu-li-sim-vra09}).

In 1972, E. Calabi \cite{cal72} found a composition formula by which one can construct new hyperbolic affine hyperspheres from any two given ones. The present author has generalized Calabi construction to the case of multiple factors (See \cite{lix93}, published in Chinese). Later in 1994 F. Dillen and L. Vrancken \cite{dil-vra94} generalized Calabi original composition to any two proper affine hyperspheres and gave a detailed study of these composed affine hyperspheres. They also mentioned that their construction applies to the case of multiple factors but with no details of it. In 2008, in order to establish their later classification mentioned above, Z.J. Hu, H.Z. Li and L. Vrancken proved one characterization of the Calabi composition of hyperbolic hyperspheres (\cite{hu-li-vra08}) in terms of special decompositions of the tangent bundle. We would like to remark that, F. Dillen, H.Z. Li and X.F. Wang has defined and studied the Calabi type composition of parallel Lagrangian submanifolds in the complex projective space $\bbc P^n$ (\cite{li-wang11}).

In this paper, we explicitly define, in a unified manner, the Calabi composition of multiple factors of hyperbolic hyperspheres, possibly including some point factors viewing as ``$0$-dimensional hyperbolic hyperspheres''. After introducing the formula of definition, we make it in detail for the computation of the basic affine invariants of this composition which turn out useful in some later application. For example, as the first application of those computation we prove that the composed affine hypersphere is symmetric if and only if each of its composing factors of positive dimension is symmetric (see Theorem \ref{added}).

{\sc Acknowledgement} The author is grateful to Professor A-M Li for his encouragement and important suggestions during this study. He also thanks Professor Z.J. Hu for providing him valuable related references some of which are listed in the end of this paper.

\section{The equiaffine geometry of hypersurfaces}

In this section, we briefly present some basic facts in the equiaffine geometry of hypersurfaces. For details the readers are referred to the text books, say, \cite{li-sim-zhao93} and \cite{nom-sas94}.

Let $x:M^n\to\bbr^{n+1}$ be nondegenerate hypersurface. Then there are several basic equiaffine invariants of $x$ among which are: the affine metric (Berwald-Blaschke metric) $g$, the affine normal $\xi:=\fr1n\Delta_gx$, the Fubini-Pick $3$-form (the so called cubic form) $A\in\bigodot^3T^*M^n$ and the affine second fundamental $2$-form $B\in\bigodot^2T^*M^n$ (See for example \cite{li-sim-zhao93} and \cite{nom-sas94}). By using the index lifting by the metric $g$, we can identify $A$ and $B$ with the linear maps $A:TM\to \ed(TM)$ or $A:TM\bigodot TM\to TM$ and $B:TM\to TM$, respectively, by
\be\label{ab}
g(A(X)Y,Z)=A(X,Y,Z) \mb{\ or\ }g(A(X,Y),Z)=A(X,Y,Z),\quad
g(B(X),Y)=B(X,Y),
\ee
for all $X,Y,Z\in TM$. Sometimes we call the corresponding $B\in \ed(TM)$ the affine shape operator of $x$. In this sense, the affine Gauss equation can be written as follows:
\be\label{gaus}
R(X,Y)Z=\fr12(g(Y,Z)B(X)+B(Y,Z)X-g(X,Z)B(Y)-B(X,Z)Y)-[A(X),A(Y)](Z),
\ee
where, for any linear transformations $T,S\in \ed(TM)$,
\be\label{comm}
[T,S]=T\circ S-S\circ T.
\ee
Each of the eigenvalues $B_1,\cdots,B_n$ of the linear map $B:TM\to TM$ is called the affine principal curvature of $x$. Define
\be\label{afme}
L_1:=\fr1n\tr B=\fr1n\sum_iB_i.
\ee
Then $L_1$ is referred to as the affine mean curvature of $x$. A hypersurface $x$ is called an (elliptic, parabolic, or hyperbolic) affine hypersphere, if all of its affine principal curvatures are equal to one (positive, 0, or negative) constant. In this case we have
\be\label{afsp}
B(X)=L_1X,\quad\mb{for all\ }X\in TM.
\ee
It follows that the affine Gauss equation \eqref{gaus} of an affine hypersphere assumes the following form:
\be\label{gaus_af sph}
R(X,Y)Z=L_1(g(Y,Z)X-g(X,Z)Y)-[A(X),A(Y)](Z),
\ee

Furthermore, all the affine lines of an elliptic affine hypersphere or a hyperbolic affine hypersphere $x:M^n\to\bbr^{n+1}$ pass through a fix point $o$ which is refer to as the affine center of $x$; Both the elliptic affine hyperspheres and the hyperbolic affine hyperspheres are called proper affine hyperspheres, while the parabolic affine hyperspheres are called improper affine hyperspheres.

{\prop\label{affine spheres} (\cite{li-sim-zhao93}) A nondegenerate hypersurface $x:M^n\to \bbr^{n+1}$ is a proper affine hypersphere with affine mean curvature $L_1$ and with the origin $o$ as its affine center if and only if the affine line is parallel to the position vector $x$. In this case, the affine normal $\xi$ is given by $\xi=-L_1x$.}

For each vector field $\eta$ transversal to the tangent space of $x$, we have the following direct decomposition
$$
x^*T\bbr^{n+1}=x_*(TM)\oplus \bbr\cdot\eta.
$$
This decomposition and the canonical differentiation $\bar D^0$ on $\bbr^{n+1}$ define a bilinear form $h\in\bigodot^2T^*M^n$ and a connection $D^\eta$ on $TM$ as follows:
\be\label{dfn h}
\bar D^0_XY=x_*(D^\eta_XY)+h(X,Y)\eta,\quad\forall X,Y\in TM.
\ee
\eqref{dfn h} can be referred as to the {\em affine Gauss formula} of the hypersurface $x$.
In particular, in case that $\eta$ is parallel to the affine normal $\xi$, the induced connection $\nabla:=D^\eta$ is independent of the choice of $\eta$ and is referred to as the affine connection of $x$.

In what follows we make the following convention for the range of indices:
$$1\leq i,j,k,l\leq n.$$

Let $\{e_i,e_{n+1}\}$ be a local unimodular frame field along $x$ with $\eta:=e_{n+1}$ parallel to the affine normal $\xi$, and $\{\omega^i,\omega^{n+1}\}$ be its dual coframe. Then we have connection forms $\omega^A_B$, $1\leq A,B\leq n+1$, defined by
$$
d\omega^A=\omega^B\wedge\omega^A_B,\quad d\omega^A_B=\sum_{C=1}^{n+1}\omega^C_A\wedge\omega^B_A,\quad \omega^{n+1}\equiv 0.
$$
Furthermore, the above invariants can be respectively expressed locally as
\be\label{gab}
g=\sum g_{ij}\omega^i\omega^j,\quad A=\sum A_{ijk}\omega^i\omega^j\omega^k,\quad B=\sum B_{ij}\omega^i\omega^j,
\ee
subject to the following basic formulas:
\begin{align}
&\sum_{i,j} g^{ij}A_{ijk}=0,\text{\ or equivalently\ }\omega^{n+1}_{n+1}+\fr1{n+2}d\log H=0,\label{basic1}\\
&A_{ijk,l}-A_{ijl,k}=\fr12(g_{ik}B_{jl}+g_{jl}B_{ik} -g_{il}B_{jk}-g_{jk}B_{il}),\label{basic3}\\
&\sum_{l}A^l_{ij,l}=\fr n2(L_1g_{ij}-B_{ij}),\label{basic3-1}
\end{align}
where $A_{ijk,l}$ are the covariant derivatives of $A_{ijk}$ with respect to Levi-Civita connection of $g$.

Write $h=\sum h_{ij}\omega^i\omega^j$ and $H=|\det(h_{ij})|$. Then
\be\label{dfn g}
g_{ij}=H^{-\fr1{n+2}}h_{ij},\quad \xi=H^{\fr1{n+2}}e_{n+1}.
\ee
Define
\be\label{hijk0}
\sum_kh_{ijk}\omega^k=dh_{ij}+h_{ij}\omega^{n+1}_{n+1}-\sum h_{kj}\omega^k_i-\sum h_{ik}\omega^k_j.
\ee
Then the Fubini-Pick form $A$ can be determined by the following formula:
\be\label{hijktoaijk}
A_{ijk}=-\fr12H^{-\fr1{n+2}}h_{ijk}.
\ee

To end this section, we would like to introduce the following concept:

{\dfn A nondegenerate hypersurface $x:M^n\to \bbr^{n+1}$ is called affine symmetric (resp. locally affine symmetric) if

(1) the pseudo-Riemannian manifold $(M^n,g)$ is symmetric (resp. locally symmetric) and therefore $(M^n,g)$ can be written (resp. locally written) as $G/K$ for some connected Lie group $G$ of isometries with $K$ one of its closed subgroups;

(2) the Fubini-Pick form $A$ is invariant under the action of $G$.}

\section{Calabi composition of multiple hyperbolic affine hyperspheres}

In this section, we aim to derive an explicit formula to define the Calabi composition of multiple factors of hyperbolic hyperspheres, possibly including $0$-dimensional factors, and then make some detailed and necessary computations for the basic affine invariants of this composition. This treatment seems not to have appeared in the literature other than \cite{lix93} where, for the first time, the author introduced the concept of multiple factor composition of the Calabi's type.

\subsection{The Calabi composition of two hyperbolic affine hyperspheres}\label{cal prod}

\newcommand{\stx}[2]{\strl{(#1)}{#2}}
\newcommand{\spec}[1]{\prod_{#1=1}^K\fr{c_{#1}^{n_{#1}+1}H_{(#1)}^{\fr1{n_{#1}+2}}} {(n_{#1}+1)(-\!\!\stx{#1}{L}_1)}}
\mb{}\par
For completeness we start with the simplest case, that is, the Calabi Composition of two factors which has appeared in many literatures.

Let $x_a:M^{n_a}\to \bbr^{n_a+1}$ be two hyperbolic affine hyperspheres with affine mean curvatures $\strl{(a)}{L}_1$, $a=1,2$,  and with the origin their common affine center.

{\prop\label{calpro} $($\cite{cal72}, \cite{li-sim-zhao93}, \cite{lix93}, \cite{dil-vra94}$)$ For any positive numbers $c_1,c_2$, define $M^n=\bbr\times M^{n_1}\times M^{n_2}$ and $x:M^n\to\bbr^{n+1}$ such that
\be\label{2factors}
x(t,p_1,p_2)=\left\{c_1\exp\left( \fr t{n_1+1}\right)x_1(p_1),c_2\exp \left(\fr {-t}{n_2+1}\right)x_2(p_2)\right\},\quad
\forall\,(t,p_1,p_2)\in M^n.
\ee
Then $x$ is again a hyperbolic affine hypersphere with affine mean curvature $L_1=-\fr1{(n+1)C}$, called {\em the Calabi composition of $x_1$ and $x_2$}, where the constant $C$ is given by
$$
C^{n+1}= \fr{c_1^{2(n_1+1)}c_2^{2(n_2+1)}} {(n_1+1)^{n_1+1}(n_2+1)^{n_2+1}(n_1+n_2+2)(-\strl{(1)}{L}_1)^{(n_1+2)} (-\strl{(2)}{L}_1)^{(n_2+2)}}.
$$

Furthermore, if $\bar x_a$ are equiaffine equivalent to $x_a$, $a=1,2$, respectively, then the Calabi composition $\bar x$ of $\bar x_1$ and $\bar x_2$ with the same constants $c_1,c_2$ is equiaffine equivalent to $x$.}

\dfn\rm The Calabi composition $x$ of $x_1$ and $x_2$ with $c_1=c_1=1$ can be referred to as {\em the Calabi product of $x_1$ and $x_2$}, and is denoted by $x=x_1*x_2$.

Then we have

{\prop\label{comml}{\rm(Almost Commutative Law)} Let $x_1,x_2$ be two hyperbolic affine hyperspheres of dimensions $n_1,n_2$ respectively. Then the Calabi products $x_2*x_1$ and $x_1*x_2$ are ``almost the same'', that is, they differ only by a linear transformation on $\bbr^{n_1+n_2+2}$ of determinant $(-1)^{(n_1+1)(n_2+1)}$. In other words, $x_2*x_1$ and $x_1*x_2$ are either equiaffine equivalent to each other, or differ by
a linear transformation on $\bbr^{n_1+n_2+2}$ of determinant $-1$.}

\proof In fact, denoting the identity matrix of order $m$ by $I_m$, we have
\begin{align}
x_2*x_1=&\left(\exp\left(\fr{t}{n_2+1}\right)x_2, \exp\left(-\fr{t}{n_1+1}\right)x_1\right)\nonumber\\
=&\left(\exp\left(\fr{t_1}{n_1+1}\right)x_1, \exp\left(-\fr{t_1}{n_2+1}\right)x_2\right)\lmx 0& \exp\left(-\fr{t+t_1}{n_1+1}\right)I_{n_1+1}\\ \exp\left(\fr{t+t_1}{n_2+1}\right)I_{n_2+1}&0\rmx\nonumber\\
=&x_1*x_2\lmx 0& \exp\left(-\fr{t+t_1}{n_1+1}\right)I_{n_1+1}\\ \exp\left(\fr{t+t_1}{n_2+1}\right)I_{n_2+1}&0\rmx\nonumber.
\end{align}
By making the coordinate change by $t_1=-t$, we find that
$$
x_2*x_1=x_1*x_2\lmx 0& I_{n_1+1}\\ I_{n_2+1}&0\rmx.
$$
Then the proposition follows since $$\det\lmx 0& I_{n_1+1}\\ I_{n_2+1}&0\rmx=(-1)^{(n_1+1)(n_2+1)}.$$

{\prop\label{assl}{\rm(Associative Law)} For any three hyperbolic affine hyperspheres $x_1,x_2,x_3$, $(x_1*x_2)*x_3$ and $x_1*(x_2*x_3)$ are equiaffine equivalent and thus we have $$(x_1*x_2)*x_3=x_1*(x_2*x_3).$$}

\proof By definition
\begin{align}
x_1*x_2=&\left(\exp\left(\fr{t_1}{n_1+1}\right)x_1, \exp\left(-\fr{t_1}{n_2+1}\right)x_2\right),\nonumber\\
x_2*x_3=&\left(\exp\left(\fr{t'_1}{n_2+1}\right)x_3, \exp\left(-\fr{t'_1}{n_3+1}\right)x_3\right).\nonumber
\end{align}
Therefore
\begin{align}
(x_1*x_2)*x_3=&\left(\exp\left(\fr{t_2}{n_1+n_2+2}\right)\cdot \exp\left(\fr{t_1}{n_1+1}\right)x_1,\right.\nonumber\\
&\left.\hs{2cm}\exp\left(\fr{t_2}{n_1+n_2+2}\right)\cdot \exp\left(\fr{-t_1}{n_2+1}\right)x_2, \exp\left(\fr{-t_2}{n_3+1}\right)x_3\right)\nonumber\\
=&\left(\exp\left(\fr{t_1}{n_1+1}+\fr{t_2}{n_1+n_2+2}\right)x_1,\right.\nonumber\\
&\left.\hs{2cm}\exp\left(-\fr{t_1}{n_2+1}+\fr{t_2}{n_1+n_2+2}\right)x_2, \exp\left(-\fr{t_2}{n_3+1}\right)x_3\right);\nonumber\\
x_1*(x_2*x_3)=&\left(\exp\left(\fr{t'_2}{n_1+1}\right)x_1, \exp\left(\fr{t'_1}{n_2+1}-\fr{t'_2}{n_2+n_3+2}\right)x_2, \right.\nonumber\\
&\left.\hs{5cm} \exp\left(-\fr{t'_1}{n_3+1}-\fr{t'_2}{n_2+n_3+2}\right)x_3\right).\nonumber
\end{align}

Now, we make the following change of coordinates:
\begin{align}
t'_1=&-\fr{n_3+1}{n_2+n_3+2}t_1
+\fr{(n_2+1)(n_1+n_2+n_3+3)}{(n_1+n_2+2)(n_2+n_3+2)}t_2,\nonumber\\
t'_2=&t_1+\fr{n_1+1}{n_1+n_2+2}t_2.\nonumber
\end{align}
Then we have
\begin{align}
&\fr{t_1}{n_1+1}+\fr{t_2}{n_1+n_2+2}=\fr{t'_2}{n_1+1},\nnm\\ &-\fr{t_1}{n_2+1}+\fr{t_2}{n_1+n_2+2}=\fr{t'_1}{n_2+1}-\fr{t'_2}{n_2+n_3+2},\nnm\\
&-\fr{t_2}{n_3+1}=-\fr{t'_1}{n_3+1}-\fr{t'_2}{n_2+n_3+2}. \nnm
\end{align}
Thus $(x_1*x_2)*x_3=x_1*(x_2*x_3)$.

\subsection{The Calabi composition of more factors---the definition and the affine metric}\label{calbi product}

\mb{}\par
We are to generalize formula \eqref{2factors}, and to this end we need some necessary notations. Firstly, given an integer $K\geq 2$ and $K$ nonnegative integers $n_1,\cdots, n_K$, we define $n=K-1+n_1+\cdots+n_K$ and make the following conventions for the ranges of different kinds of indices:
\begin{align}
&1\leq a,b,c\cdots\leq K,\quad 1\leq\lambda,\mu,\nu\leq K-1,\quad 1\leq i_a,j_a,k_a\leq n_a,\\
&\bar i_a=i_a+K-1+\sum_{b<a}n_b,\  \bar j_a=j_a+K-1+\sum_{b<a}n_b,\ \bar k_a=k_a+K-1+\sum_{b<a}n_b.
\end{align}

Secondly, for each $a=1,2,\cdots,K$ and $(t_1,\cdots,t_{K-1})\in\bbr^{K-1}$, we define $f_a: =n_1+\cdots+n_a+a$ and
$$e_a:=\exp\left(-\fr{t_{a-1}}{n_{a}+1}+\fr{t_{a}}{f_{a}}+\fr{t_{a+1}}{f_{a+1}} +\cdots+\fr{t_{K-1}}{f_{K-1}}\right).$$
In particular,
$$
e_1=\exp\left(\fr{t_1}{f_1}+\fr{t_2}{f_2} +\cdots+\fr{t_{K-1}}{f_{K-1}}\right),\quad
e_K=\exp\left(-\fr{t_{K-1}}{n_K+1}\right).
$$

\newcommand{\la}{\stx{a}{L}\!\!_1{}}\newcommand{\lb}{\stx{b}{L}\!\!_1{}}
\newcommand{\lc}{\stx{c}{L}\!\!_1{}}\newcommand{\lalp}{\stx{\alpha}{L}\!\!_1{}}
\newcommand{\ha}{\!\stx{a}{h}{}\!\!} %\newcommand{\hb}{\stx{b}{h}\!\!{}}
\newcommand{\Ha}{H_{(a)}}\newcommand{\Hb}{H_{(b)}}\newcommand{\Hc}{H_{(c)}}
\newcommand{\ga}{\!\!\stx{a}{g}{}\!\!\!}\newcommand{\Ga}{\stx{a}{G}\!\!{}}
\newcommand{\gb}{\!\!\stx{b}{g}{}\!\!\!}\newcommand{\Gb}{\stx{b}{G}\!\!{}}
\newcommand{\galp}{\!\!\stx{\alpha}{g}{}\!\!\!}
Then our first generalization is as follows:

{\thm\label{multicomp} Let $x_a:M^{n_a}\to\bbr^{n_a+1}$, $a=1,2,\cdots,K$, be $K(\geq 2)$ hyperbolic affine hyperspheres with affine mean curvatures $\la$ and with the origin their common affine center. Then for any $K$ positive numbers $c_a$, $a=1,\cdots, K$, we have a new hyperbolic affine hypersphere $x:M^n\to\bbr^{n+1}$ with the affine mean curvature
\be\label{l1c}
L_1=\fr1{f_KC},\quad C:=\left(f^{-1}_K \prod_{a=1}^K\fr{c_a^{2(n_a+1)}} {(n_a+1)^{n_a+1}(-\!\!\stx{a}{L}_1)^{n_a+2}}\right)^{\fr1{n+2}},
\ee
where $n=f_K-1$, $M^n=\bbr^{K-1}\times M^{n_1}\times\cdots\times M^{n_K}$ and
\begin{align}
x(t^1,&\cdots,t^{K-1},p_1,\cdots,p_K):=(c_1e_1x_1(p_1), c_2e_2x_2(p_2),\cdots,c_Ke_Kx_K(p_K)),\nnm\\&\hs{1cm}\forall (t^1,\cdots,t^{K-1},p_1,\cdots,p_K)\in M^n.\label{mulpro1}
\end{align}

Moreover, for given positive numbers $c_1,\cdots,c_K$, there exits some $c>0$ and $c'>0$ such that
the following three hyperbolic affine hyperspheres
\begin{align} &x:=(c_1e_1x_1, c_2e_2x_2,\cdots,c_Ke_Kx_K),\nnm\\
&\bar x:=c(e_1x_1, e_2x_2,\cdots,e_Kx_K),\nnm\\
&\td x:=(e_1x_1, e_2x_2,\cdots,c'e_Kx_K)\nnm
\end{align}
are equiaffine equivalent to each other.}

\dfn\rm The hyperbolic affine hypersphere $x$ in the above proposition is called the Calabi composition of the given hyperbolic affine hyperspheres $x_a$, $a=1,\cdots,K$.

{\it Proof of Theorem \ref{multicomp}}

Let $(v^{i_a},\ 1\leq i_a\leq n_a)$ be a local coordinate system on $M^{n_a}$.
\newcommand{\xai}{x_{a,i_a}} \newcommand{\xaij}{x_{a,i_aj_a}}
Define
$$
u^\lambda=t^\lambda,\quad u^{\bar i_a}=v^{i_a};\quad \xai=\pp{x_a}{v^{i_a}},\quad \xaij=\ppp{x_a}{v^{i_a}}{v^{j_a}};\quad\ha_{i_aj_a}=\det(x_{a,1},\cdots,x_{a,n_a},\xaij).
$$
Then the affine metric $\ga=\sum_{i_a,j_a}\ga_{i_aj_a}dv^{i_a}dv^{j_a}$ is given by
$$
\ga_{i_aj_a}=\Ha^{-\fr1{n_a+2}}\ha_{i_aj_a},\quad \Ha:=|\det(\ha_{i_aj_a})|.
$$
If we denote $h_{ij}=\det(\pp{x}{u^1},\cdots,\pp{x}{u^n},\ppp{x}{u^i}{u^j})$, then a long but direct computation, using Proposition \ref{affine spheres} and the definition of the affine normals $\xi_a$ ($1\leq a\leq K$), concludes that
\begin{align}
h_{\lambda\mu}=&\veps\spec{a}\fr{f_{\lambda+1}} {(n_{\lambda+1}+1)f_\lambda}\delta_{\lambda\mu},
\quad h_{\lambda \bar i_a}=0,\label{h lambda mu}\\
h_{\bar i_a\bar j_b}=&\veps\spec{c}\fr{(n_a+1)(-\!\!\la)}{\Ha^{\fr1{n_a+2}}}\ha_{i_aj_a}\delta_{ab}, \label{h ia jb} \end{align}
where $\veps=(-1)^{(K+1)n_1+Kn_2+\cdots+3n_{K-1}}$. By suitably choosing the orientation of $\bbr^{K-1}$, if necessary, we can reasonably take $\veps=1$. It then follows that
\begin{align}
H:=&|\det(h_{ij})_{1\leq i,j\leq n}|=\left|\det\,{\rm diag}\left(h_{11},\cdots,h_{K-1\,K-1}, (\stx{1}{h}\!\!_{i_1j_1})_{n_1\times n_1},\cdots,(\stx{K}{h}\!\!_{i_Kj_K})_{n_K\times n_K}\right)\right|\nnm\\
=&\prod_\lambda\left(\spec{a}\fr{f_{\lambda+1}}{(n_{\lambda+1}+1)f_\lambda}\right) \cdot\prod_a\left(\spec{b}\fr{(n_a+1)(-\!\!\la)}{\Ha^{\fr1{n_a+2}}}\right)^{n_a} |\det(\ha_{i_aj_a})|\nnm\\
=&\left(\spec{a}\right)^{K-1}\cdot \prod_\lambda\fr{f_{\lambda+1}}{(n_{\lambda+1}+1)f_\lambda}\times \nnm\\
&\hs{4cm}\times\left(\spec{b}\right)^{\sum_an_a}\cdot \prod_a\left(\fr{(n_a+1)(-\!\!\la)}{\Ha^{\fr1{n_a+2}}}\right)^{n_a} \Ha\nnm\\
=&\left(\spec{a}\right)^{\sum_an_a+K-1}\cdot \prod_\lambda\fr{f_{\lambda+1}}{(n_{\lambda+1}+1)f_\lambda} \cdot\prod_a\left(\fr{(n_a+1)(-\!\!\la)}{\Ha^{\fr1{n_a+2}}}\right)^{n_a} \Ha\nnm\\
=&\prod_a\fr{c_a^{(n_a+1)(f_K-1)}\Ha^{\fr{f_K+1}{n_a+2}}} {(n_a+1)^{f_K-n_a-1}(-\!\!\la)^{f_K-n_a-1}} \cdot\prod_\lambda\fr{f_{\lambda+1}}{(n_{\lambda+1}+1)f_\lambda}
\end{align}
Since $f_1=n_1+1$, we have
\be\label{fml1}
\prod_\lambda\fr{f_{\lambda+1}}{(n_{\lambda+1}+1)f_\lambda} =\fr{f_2}{(n_2+1)f_1}\cdot\fr{f_3}{(n_3+1)f_2}\cdot\cdots\cdots\fr{f_K}{(n_K+1)f_{K-1}}
=f_K\left(\prod_a(n_a+1)\right)^{-1}.
\ee
Note that $f_K=n+1$. It follows that
\be\label{H}
H=f_K\prod_a\fr{c_a^{(n_a+1)(f_K-1)}\Ha^{\fr{f_K+1}{n_a+2}}} {(n_a+1)^{f_K-n_a}(-\!\!\la)^{f_K-n_a-1}},\quad
H^{-\fr1{n+2}}=f^{-\fr1{n+2}}_K \prod_a\fr {(n_a+1)^{\fr{f_K-n_a}{f_K+1}}(-\!\!\la)^{\fr{f_K-n_a-1}{f_K+1}}} {\left(c_a^{n_a+1}\right)^{\fr{f_K-1}{f_K+1}}\Ha^{\fr1{n_a+2}}}.
\ee
\newcommand{\HH}{f_K\prod_a\fr{c_a^{(n_a+1)(f_K-1)}\Ha^{\fr{f_K+1}{n_a+2}}}
{(n_a+1)^{f_K-n_a}(-\!\!\la)^{f_K-n_a-1}}}
\newcommand{\h}{f^{-\fr1{n+2}}_K \prod_a\fr {(n_a+1)^{\fr{f_K-n_a}{f_K+1}}(-\!\!\la)^{\fr{f_K-n_a-1}{f_K+1}}}
{\left(c_a^{n_a+1}\right)^{\fr{f_K-1}{f_K+1}}\Ha^{\fr1{n_a+2}}}}
Therefore the affine metric $g=\sum_{ij}g_{ij}$ of the hypersurface $x$ is given by
\begin{align}
g_{\lambda\mu}=&H^{-\fr1{n+2}}h_{\lambda\mu}\nnm\\
=&f^{-\fr1{n+2}}_K \prod_a\fr {(n_a+1)^{\fr{f_K-n_a}{f_K+1}}(-\!\!\la)^{\fr{f_K-n_a-1}{f_K+1}}} {\left(c_a^{n_a+1}\right)^{\fr{f_K-1}{f_K+1}}\Ha^{\fr1{n_a+2}}} \cdot \spec{a}\fr{f_{\lambda+1}}{(n_{\lambda+1}+1)f_\lambda}\delta_{\lambda\mu}\nnm\\ =&\fr{f_{\lambda+1}C}{(n_{\lambda}+1)f_\lambda}\delta_{\lambda\mu},\label{g-lmdmu}\\
g_{\bar i_a\bar j_b}=&H^{-\fr1{n+2}}h_{i_aj_b}\nnm\\
=&f^{-\fr1{n+2}}_K \prod_b\fr {(n_b+1)^{\fr{f_K-n_b}{f_K+1}}(-\!\!\stx{b}{L}\!\!_1)^{\fr{f_K-n_b-1}{f_K+1}}} {\left(c_b^{n_b+1}\right)^{\fr{f_K-1}{f_K+1}}H_{(b)}^{\fr1{n_b+2}}} \cdot \spec{b}\fr{(n_a+1)(-\!\!\la)}{\Ha^{\fr1{n_a+2}}}\ha_{i_aj_a}\delta_{ab}\nnm\\ =&(n_a+1)(-\!\!\la)C\ \ga_{i_aj_a}\delta_{ab},\label{g-iajb},\\
g_{\lambda\bar i_a}=&0,\label{g-lmdia}
\end{align}
in which the constant $C$ is defined by \eqref{l1c}. It follows that
\begin{align}
G:=&\det(g_{ij})=\det{\rm diag}(g_{11},\cdots,g_{K-1\,K-1},(g_{\bar i_1\bar j_1}),\cdots (g_{\bar i_K\bar j_K}))
=g_{11}\cdots g_{K-1\,K-1}\prod_a\det(g_{\bar i_a\bar j_a})\nnm\\ =&C^{K-1}\prod_\lambda\fr{f_{\lambda+1}}{(n_{\lambda+1}+1)f_\lambda} \cdot\prod_a\left(\left((n_a+1)(-\!\!\la)\right)^{n_a}\cdot C^{n_a}\det(\ga_{i_aj_a})\right)\nnm\\
=&C^{f_K-1}\prod_\lambda\fr{f_{\lambda+1}}{(n_{\lambda+1}+1)f_\lambda} \cdot\prod_a\left(\left((n_a+1)(-\!\!\la)\right)^{n_a}\Ga\right)\nnm\\
=&C^{f_K-1}f_K\prod_a\left((n_a+1)^{n_a-1}(-\!\!\la)^{n_a}\Ga\right),\nnm
\end{align}
where the last equality uses \eqref{fml1}.

On the other hand, from \eqref{g-lmdmu}, \eqref{g-iajb} and \eqref{g-lmdia} we find that
\be\label{co-g}
g^{\lambda\mu}=\fr{(n_{\lambda}+1)f_\lambda}{f_{\lambda+1}C}\delta^{\lambda\mu}, \quad g^{\bar i_a\bar j_b}=\left((n_a+1)(-\!\!\la)C\right)^{-1}\ga\,^{i_aj_a}\delta^{ab},\quad g^{\lambda\bar i_a}=0.
\ee
Note that
\be\label{fml2}
\pp{e_a}{t^\lambda}=e_a\cdot\begin{cases} 0,&{\rm if}\ 0\leq\lambda\leq a-2;\\
-\fr1{n_a+1},&{\rm if}\ \lambda=a-1;\\
\fr1{f_\lambda},&{\rm if}\ a\leq\lambda\leq K-1.
\end{cases}
\ee
Thus
\be\label{fml3}
\fr{\partial^2e_a}{(\partial t^\lambda)^2}
=e_a\cdot\begin{cases} 0,&{\rm if}\ 0\leq\lambda\leq a-2;\\
\fr1{(n_a+1)^2},&{\rm if}\ \lambda=a-1;\\
\fr1{f^2_\lambda},&{\rm if}\ a\leq\lambda\leq K-1.
\end{cases}
\ee
Consequently
\begin{align}
\sum_\lambda g^{\lambda\lambda}\fr{\partial^2e_a}{(\partial t^\lambda)^2} =&g^{a-1\,a-1}\cdot\fr1{(n_a+1)^2}e_a +\sum_{\lambda\geq a}g^{\lambda\lambda}\fr1{f^2_\lambda}e_a\nnm\\
=&\fr{e_a}{C}\left(\fr{f_{a-1}}{(n_a+1)f_a} +\sum_{\lambda\geq a}\fr{n_{\lambda+1}+1}{f_{\lambda+1}f_\lambda}\right)\nnm\\
=&\fr{e_a}{C}\left(\fr{f_{a-1}}{(n_a+1)f_a} +\sum_{\lambda\geq a}\left(\fr1{f_\lambda}-\fr1{f_{\lambda+1}}\right)\right)\nnm\\
=&\fr{e_a}{C}\left(\left(\fr1{n_a+1}-\fr1{f_a}\right) +\left(\fr1{f_a}-\fr1{f_K}\right)\right)\nnm\\
=&\fr{e_a}{C}\left(\fr1{n_a+1}-\fr1{f_K}\right).\label{fml3'}
\end{align}
Moreover, for each $a$
\begin{align}
&\fr1{\sqrt{|G|}}\sum_{i_a,j_a}\pp{}{u^{\bar i_a}}\left(g^{\bar i_a\bar j_a}\sqrt{|G|}c_ae_a
\pp{x_a}{u^{\bar j_a}}\right)
=c_ae_a\left(C^{f_K-1}f_K\prod_b\left((n_b+1)^{n_b-1} (-\!\!\lb)^{n_b}|\Gb|\right)\right)^{-\fr12}\times\nnm\\
&\hs{.5cm}\times \sum_{i_a,j_a}\pp{}{v^{i_a}} \left(\left(C^{f_K-1}f_K\prod_b\left((n_b+1)^{n_b-1} (-\!\!\lb)^{n_b}|\Gb|\right)\right)^{\fr12}\left((n_a+1) (-\!\!\la)C\right)^{-1}\ga\, ^{i_aj_a}\pp{x_a}{v^{j_a}}\right)\nnm\\
=&\fr{c_ae_a}{(n_a+1) (-\!\!\la)C}\fr1{\sqrt{|\Ga|}}\sum_{i_a,j_a}\pp{}{v^{i_a}}\left(\ga\, ^{i_aj_a}\sqrt{|\Ga|}\pp{x_a}{v^{j_a}}\right)\nnm\\
=&\fr{c_ae_a}{(n_a+1) (-\!\!\la)C}\Delta\ _{\ga}\ (x_a)=\fr{n_a}{(n_a+1)C}c_ae_ax_a,\label{fml4}
\end{align}
where we have used the fact that (see Proposition \ref{affine spheres})
\be\label{fml5}
\Delta\ _{\ga}\ (x_a)=n_a\xi_a={n_a}(-\!\!\la)x_a.
\ee

Now we can use \eqref{fml3'} and \eqref{fml4} to compute the Laplacian of $x$ with the metric $g$ as follows:
\begin{align}
\Delta_g\,(x)=&\sum_{i,j}\fr1{\sqrt{|G|}}\pp{}{u^i} \left(\sqrt{|G|}g^{ij}\pp{x}{u^j}\right)\nnm\\
=&\fr1{\sqrt{|G|}}\sum_{\lambda}\pp{}{t^\lambda} \left(\sqrt{|G|}g^{\lambda\lambda}\pp{x}{t^\lambda}\right) +\sum_{a,i_a,j_a}\fr1{\sqrt{|G|}}\pp{}{u^{\bar i_a}}\left(\sqrt{|G|}g^{\bar i_a\bar j_a}\pp{x}{u^{\bar j_a}}\right)\nnm\\
=&\left(\cdots,c_ax_a\sum_\lambda\fr1{\sqrt{|G|}}\pp{}{t^\lambda} \left(\sqrt{|G|}g^{\lambda\lambda}\pp{e_a}{t^\lambda}\right) +\fr1{\sqrt{|G|}}\pp{}{u^{\bar i_a}}\left(g^{\bar i_a\bar j_a}\sqrt{|G|}c_ae_a
\pp{x_a}{u^{\bar j_a}}\right),\cdots\right)\nnm\\
=&\left(\cdots,c_ax_a\fr{e_a}{C}\left(\fr1{n_a+1}-\fr1{f_K}\right) +\fr{n_a}{(n_a+1)C}c_ae_ax_a,\cdots\right)\nnm\\
=&\left(\cdots,-\fr{f_{K-1}}{f_KC}c_ae_ax_a,\cdots\right)=-\fr{f_{K-1}}{f_KC}x.
\end{align}
Thus we find the affine normal
\be\label{afnml}
\xi=\fr1{n}\Delta_g\,x=-\fr1{f_KC}x,
\ee
which, via Propostion \ref{affine spheres}, proves that $x$ is a hyperbolic affine hypersphere with affine mean curvature
\be\label{afmn}
L_1=-\fr1{f_KC}.
\ee

Finally, for any set of $c_1,\cdots,c_K$, take $c=\left(\prod_ac_a^{n_a+1}\right)^{\fr1{f_K}}$ and $c'=\prod_ac_a^{n_a+1}$. Then it is easily seen that
$x$, $\bar x$ and $\td x$ are equiaffine equivalent each other.\endproof

\subsection{The induced affine connection}

\mb{}\par
Later in this section, we shall compute the Fubini-Pick form of $x$. To this end, we first need to find the induced affine connection or, equivalently, the corresponding connection forms $\omega^j_i=\Gamma^j_{ik}\omega^k$, where $\{\omega^i,\omega^{n+1}\}$ is the dual coframe of the unimodular frame $\{\pp{x}{u^i},e_{n+1}\}$.

Since $x_a$ and $x$ are hyperbolic affine hyperspheres, we have
\newcommand{\gma}{\stx{a}{\Gamma}{}\!\!}\newcommand{\gmb}{\stx{b}{\Gamma}{}}
\be\label{fml6}
\ppp{x_a}{v^{i_a}}{v^{j_a}}+\ga_{i_aj_a}\!\la x_a=\sum_{k_a}\gma^{k_a}_{i_aj_a}\pp{x_a}{v^{k_a}},\quad \ppp{x}{u^i}{u^j}+g_{ij}L_1 x=\sum_{k}\Gamma^k_{ij}\pp{x}{u^k}.
\ee

{\lem\label{gamij} It holds that
\be\label{gmij}
\Gamma^{\bar k_b}_{\bar i_a\bar j_a}=\gma^{k_a}_{i_aj_a}\delta^b_a,\quad
\Gamma^{\lambda}_{\bar i_a\bar j_a}=\begin{cases} 0,&{\rm if\ }1\leq\lambda\leq a-2;\\
\fr{(n_a+1)f_{a-1}}{f_a}\!\la\ga_{i_aj_a},&{\rm if\ }\lambda=a-1;\\
-\fr{(n_a+1)(n_{\lambda+1}+1)}{f_{\lambda+1}}\!\la\ga_{i_aj_a},&{\rm if\ }a\leq\lambda\leq K-1.
\end{cases}
\ee}

\proof Directly we compute
\begin{align}
&\ppp{x}{u^{\bar i_a}}{u^{\bar j_a}}+g_{\bar i_a\bar j_a}L_1 x\nnm\\ =&\left(0,\cdots,c_ae_a\ppp{x_a}{v^{i_a}}{v^{j_a}},\cdots,0\right) +(n_a+1)\big(-\!\!\la\big)C\ga_{i_aj_a}\left(-\fr1{f_KC}\right) (\cdots,c_be_bx_b,\cdots)\nnm\\
=&\left(0,\cdots,c_ae_a\ga_{i_aj_a}\!\la x_a +\sum_{k_a}\gma^{k_a}_{i_aj_a}\pp{x_a}{v^{k_a}}, \cdots,0\right) +\fr{(n_a+1)\!\la}{f_K}\ga_{i_aj_a}(\cdots,c_be_bx_b,\cdots)\nnm\\
=&\left(c_1e_1\fr{(n_a+1)\!\la}{f_K}\ga_{i_aj_a}x_1,\cdots, c_ae_a\!\la \left(\fr{(n_a+1)\!\la}{f_K}-1\right) \ga_{i_aj_a}x_a+c_ae_a\sum_{k_a}\gma^{k_a}_{i_aj_a}\pp{x_a}{v^{k_a}},\right.\nnm\\
&\hs{8cm}\left.\cdots,c_Ke_K\fr{(n_a+1) \!\la}{f_K}\ga_{i_aj_a}x_K\right).\label{fml8}
\end{align}
On the other hand, by using \eqref{fml2} it is easy to find that, for each $b=1,2,\cdots,K$
\begin{align}
\sum_k\Gamma^k_{\bar i_a\bar j_a}\pp{x}{u^k} =&\sum_\lambda\Gamma^\lambda_{\bar i_a\bar j_a}\pp{x}{t^\lambda}+\sum_{b,k_b}\Gamma^{\bar k_b}_{\bar i_a\bar j_a}\pp{x}{u^{\bar k_b}}\nnm\\
=&\left(\cdots,c_b\sum_\lambda\Gamma^\lambda_{\bar i_a\bar j_a}\pp{e_b}{t^\lambda}x_b,\cdots\right) +\sum_{b,k_b}\Gamma^{\bar k_b}_{\bar i_a\bar j_b}\left(0,\cdots,c_be_b\pp{x_b}{v^{k_b}},\cdots\right)\nnm\\
=&\left(\cdots,\left(-\fr1{n_b+1}\Gamma^{b-1}_{\bar i_a\bar j_a}+\sum_{\lambda\geq b}\fr1{f_\lambda}\Gamma^\lambda_{\bar i_a\bar j_a}\right)c_be_bx_b,\cdots\right) +\left(\cdots,c_be_b\sum_{k_b}\Gamma^{\bar k_b}_{\bar i_a\bar j_a}\pp{x_b}{v^{k_b}},\cdots\right)\nnm\\
=&\left(\cdots,\left(-\fr1{n_b+1}\Gamma^{b-1}_{\bar i_a\bar j_a}+\sum_{\lambda\geq b}\fr1{f_\lambda}\Gamma^\lambda_{\bar i_a\bar j_a}\right)c_be_bx_b+c_be_b\sum_{k_b}\Gamma^{\bar k_b}_{\bar i_a\bar j_a}\pp{x_b}{v^{k_b}},\cdots\right).\label{fml9-0}
\end{align}
It then follows from \eqref{fml6}, \eqref{fml8} and \eqref{fml9-0} that
\begin{align}
&\left(-\fr1{n_b+1}\Gamma^{b-1}_{\bar i_a\bar j_a}+\sum_{\lambda\geq b}\fr1{f_\lambda}\Gamma^\lambda_{\bar i_a\bar j_a}\right)c_be_bx_b+c_be_b\sum_{k_b}\Gamma^{\bar k_b}_{\bar i_a\bar j_a}\pp{x_b}{v^{k_b}}\nnm\\
=&\begin{cases} c_be_b\fr{(n_a+1)\,\la}{f_b}\,\ga_{i_aj_a}x_b,&{\rm if\ } b\neq a;\\
c_ae_a\!\la \left(\fr{n_a+1}{f_K}-1\right) \ga_{i_aj_a}x_a+c_ae_a\sum_{k_a}\gma^{k_a}_{i_aj_a}\pp{x_a}{v^{k_a}}, &{\rm if\ } b=a
\end{cases}\label{gmij0}
\end{align}
From \eqref{gmij0} we find that
$$
\Gamma^{\bar k_a}_{\bar i_a\bar j_a}=\gma^{k_a}_{i_aj_a}, \quad \Gamma^{\bar k_b}_{\bar i_a\bar j_a}=0,\ {\rm if\ }b\neq a.
$$
Moreover, \eqref{gmij0} also implies the following linear system for $\Gamma^\lambda_{\bar i_a\bar j_a}$:
\begin{align}
&\lmx
\fr1{n_1+1}&\fr1{f_2}&\fr1{f_3}&\cdots&\fr1{f_{a-1}}&\fr1{f_a}&\fr1{f_{a+1}}& \cdots&\fr1{f_{K-2}}&\fr1{f_{K-1}}\\
-\fr1{n_2+1}&\fr1{f_2}&\fr1{f_3}&\cdots&\fr1{f_{a-1}}&\fr1{f_a}&\fr1{f_{a+1}}& \cdots&\fr1{f_{K-2}}&\fr1{f_{K-1}}\\
0&-\fr1{n_3+1}&\fr1{f_3}&\cdots&\fr1{f_{a-1}}&\fr1{f_a}&\fr1{f_{a+1}}& \cdots&\fr1{f_{K-2}}&\fr1{f_{K-1}}\\
\cdots&\cdots&\cdots&\cdots&\cdots&\cdots&\cdots&\cdots&\cdots&\cdots\\
0&0&0&\cdots&\fr1{f_{a-1}}&\fr1{f_a}&\fr1{f_{a+1}}& \cdots&\fr1{f_{K-2}}&\fr1{f_{K-1}}\\
0&0&0&\cdots&-\fr1{n_a+1}&\fr1{f_a}&\fr1{f_{a+1}}& \cdots&\fr1{f_{K-2}}&\fr1{f_{K-1}}\\
0&0&0&\cdots&0&-\fr1{n_{a+1}+1}&\fr1{f_{a+1}}& \cdots&\fr1{f_{K-2}}&\fr1{f_{K-1}}\\
\cdots&\cdots&\cdots&\cdots&\cdots&\cdots&\cdots&\cdots&\cdots&\cdots\\
0&0&0&\cdots&0&0&0& \cdots&-\fr1{n_{K-1}+1}&\fr1{f_{K-1}}\\
0&0&0&\cdots&0&0&0& \cdots&0&-\fr1{n_K+1}
\rmx
\lmx
\Gamma^1_{\bar i_a\bar j_a}\\ \Gamma^2_{\bar i_a\bar j_a}\\ \Gamma^3_{\bar i_a\bar j_a}\\
 \\ \vdots \\ \Gamma^\lambda_{\bar i_a\bar j_a}\\ \vdots \\ \\ \Gamma^{K-1}_{\bar i_a\bar j_a}
\rmx
\nnm\\
=&\left(\fr{n_a+1}{f_K}\,\la\ga_{i_aj_a}, \fr{n_a+1}{f_K} \,\la\ga_{i_aj_a},\cdots,\right.\nnm\\ &\hs{1cm}\left.\fr{n_a+1}{f_K}\,\la\ga_{i_aj_a}, \left(\fr{n_a+1}{f_K}-1\right) \la\ga_{i_aj_a},\cdots, \fr{n_a+1}{f_K} \,\la\ga_{i_aj_a}\right)^t.\label{gmij1}
\end{align}
It is not hard to show that the linear system \eqref{gmij1} is equivalent to
\begin{align}
&\lmx
\fr{f_2}{(n_1+1)(n_2+1)}\!\!\!\!\!\!&0\!\!\!\!\!\!&0\!\!\!\!\!\!&\cdots \!\!\!\!\!\!&0\!\!\!\!\!\!&0\!\!\!\!\!\!&0\!\!\!\!\!\!& \cdots\!\!\!\!\!\!&0\!\!\!\!&0\\
-\fr1{n_2+1}\!\!\!\!\!\!&\fr{f_3}{(n_3+1)f_2}\!\!\!\!\!\!&0\!\!\!\!\!\!&\cdots \!\!\!\!\!\!&0\!\!\!\!\!\!&0\!\!\!\!\!\!&0\!\!\!\!\!\!& \cdots\!\!\!\!\!\!&0\!\!\!\!&0\\
0\!\!\!\!\!\!&-\fr1{n_3+1}\!\!\!\!\!\!&\fr{f_4}{(n_4+1)f_3}\!\!\!\!\!\!&\cdots \!\!\!\!\!\!&0\!\!\!\!\!\!&0\!\!\!\!\!\!&0\!\!\!\!\!\!& \cdots\!\!\!\!\!\!&0\!\!\!\!&0\\
\cdots\!\!\!\!\!\!&\cdots\!\!\!\!\!\!&\cdots\!\!\!\!\!\!&\cdots\!\!\!\!\!\!&\cdots \!\!\!\!\!\!&\cdots\!\!\!\!\!\!&\cdots\!\!\!\!\!\!&\cdots\!\!\!\!\!\!&\cdots\!\!\!\!&\cdots\\
0\!\!\!\!\!\!&0\!\!\!\!\!\!&0\!\!\!\!\!\!&\cdots \!\!\!\!\!\!&\fr{f_a}{(n_a+1)f_{a-1}}\!\!\!\!\!\!&0\!\!\!\!\!\!&0\!\!\!\!\!\!& \cdots\!\!\!\!\!\!&0\!\!\!\!&0\\
0\!\!\!\!\!\!&0\!\!\!\!\!\!&0\!\!\!\!\!\!&\cdots \!\!\!\!\!\!&-\fr1{n_a+1}\!\!\!\!\!\!&\fr{f_{a+1}}{(n_{a+1}+1)f_a}\!\!\!\!\!\!&0\!\!\!\!\!\!& \cdots\!\!\!\!\!\!&0\!\!\!\!&0\\
0\!\!\!\!\!\!&0\!\!\!\!\!\!&0\!\!\!\!\!\!&\cdots \!\!\!\!\!\!&0\!\!\!\!\!\!&-\fr1{n_{a+1}+1}\!\!\!\!\!\!&\fr{f_{a+2}}{(n_{a+2}+1)f_{a+1}}\!\!\!\!\!\!& \cdots\!\!\!\!\!\!&0\!\!\!\!&0\\
\cdots\!\!\!\!\!\!&\cdots\!\!\!\!\!\!&\cdots\!\!\!\!\!\!&\cdots \!\!\!\!\!\!&\cdots\!\!\!\!\!\!&\cdots\!\!\!\!\!\!&\cdots\!\!\!\!\!\!&\cdots \!\!\!\!\!\!&\cdots\!\!\!\!&\cdots\\
0\!\!\!\!\!\!&0\!\!\!\!\!\!&0\!\!\!\!\!\!&\cdots\!\!\!\!\!\!&0\!\!\!\!\!\!&0\!\!\!\!\!\!&0\!\!\!\!\!\!& \cdots\!\!\!\!\!\!&-\fr1{n_{K-1}+1}\!\!\!\!&\fr{f_K}{(n_K+1)f_{K-1}}\\
0\!\!\!\!\!\!&0\!\!\!\!\!\!&0\!\!\!\!\!\!&\cdots\!\!\!\!\!\!&0\!\!\!\!\!\!&0\!\!\!\!\!\!&0\!\!\!\!\!\!& \cdots\!\!\!\!\!\!&0\!\!\!\!&-\fr1{n_K+1}
\rmx
\lmx
\Gamma^1_{\bar i_a\bar j_a}\\ \Gamma^2_{\bar i_a\bar j_a}\\ \Gamma^3_{\bar i_a\bar j_a}\\
\vdots\\ \vdots\\ \Gamma^\lambda_{\bar i_a\bar j_a}\\ \vdots\\ \vdots\\ \Gamma^{K-1}_{\bar i_a\bar j_a}
\rmx
\nnm\\
=&\left(0,0,\cdots,\,\la\ga_{i_aj_a},-\la\ga_{i_aj_a}, 0,\cdots,0,\fr{n_a+1}{f_K}\la\ga_{i_aj_a}\right)^t.\label{gmij2}
\end{align}
From \eqref{gmij2} we easily find each $\Gamma^\lambda_{i_aj_a}$ as listed in \eqref{gmij}.
\endproof

{\lem\label{gamlammu} It holds that
\begin{align}
&\Gamma^{\nu}_{\lambda\mu}=\delta^\nu_\lambda\fr1{f_\mu},\quad{\rm for\ }\lambda<\mu;\quad
\Gamma^{\nu}_{\lambda\mu}=\delta^\nu_\mu\fr1{f_\lambda},\quad{\rm for\ }\lambda>\mu,\label{gmlmdmu}\\
&\Gamma^\lambda_{\lambda\lambda}=\fr1{f_\lambda}-\fr1{n_{\lambda+1}+1},\\
&\Gamma^{\mu}_{\lambda\lambda}=\begin{cases} \fr{(n_{\mu+1}+1)f_{\lambda+1}} {(n_{\lambda+1}+1)f_{\mu+1}f_\lambda},&{\rm for\ }\lambda+1\leq\mu\leq K-1,\\
0,& {\rm otherwise};\end{cases}\\
&\Gamma^{\bar i_a}_{\lambda\mu}=0.\label{gmlmdmu1}
\end{align}}

\proof Case (1) For $\lambda<\mu$, Consider
$$
\ppp{x}{t^\lambda}{t^\mu}+g_{\lambda\mu}L_1x =\sum_\nu\Gamma^\nu_{\lambda\mu}\pp{x}{t^\nu}+\sum_{i_a}\Gamma^{\bar i_a}_{\lambda\mu}\pp{x}{u^{\bar i_a}}.
$$
Similar to the case in the proof of Lemma \ref{gamij}, we can find that $\Gamma^{\bar i_a}_{\lambda\mu}=0$ and, for $\lambda<\mu$,
\be\label{gmlmdmu0}
-\fr1{n_a+1}\Gamma^{a-1}_{\lambda\mu}+\sum_{\nu=a}^{K-1}\fr1{f_\nu}\Gamma^\nu_{\lambda\mu} =\begin{cases}\fr1{f_\lambda f_\mu},&{\rm if\ }1\leq a\leq \lambda;\\
-\fr1{(n_{\lambda+1}+1)f_\mu},&{\rm if\ }a=\lambda+1;\\
0,&{\rm if\ }a\geq\lambda+2.
\end{cases}
\ee
It follows that $\Gamma^\nu_{\lambda\mu}$, $\nu=1,\cdots,K-1$, satisfy a linear system same as \eqref{gmij1} with the right hand side being replaced by
$$
\left(\fr1{f_\lambda f_\mu},\cdots,\fr1{f_\lambda f_\mu}, -\fr1{(n_{\lambda+1}+1)f_\mu},0,\cdots,0\right)^t.
$$
This new system is equivalent to the system \eqref{gmij2} with the right hand side being replaced by
$$
\left(0,\cdots,0, \fr{f_{\lambda+1}}{(n_{\lambda+1}+1)f_\lambda f_\mu},-\fr1{(n_{\lambda+1}+1)f_\mu},\cdots,0\right)^t.
$$
Thus we easily find the expression of each $\Gamma^\nu_{\lambda\mu}$ as in \eqref{gmlmdmu}.

Case (2) For $\mu=\lambda$, Consider
$$
\fr{\partial^2x}{(\partial t^\lambda)^2}+g_{\lambda\lambda}L_1x =\sum_\mu\Gamma^\mu_{\lambda\lambda}\pp{x}{t^\mu}
+\sum_{i_a}\Gamma^{\bar i_a}_{\lambda\lambda}\pp{x}{u^{\bar i_a}},
$$
which shows in the same way that $\Gamma^{\bar i_a}_{\lambda\lambda}=0$ and
\be\label{gmlmd0}
-\fr1{n_a+1}\Gamma^{a-1}_{\lambda\lambda} +\sum_{\mu=a}^{K-1}\fr1{f_\mu}\Gamma^\mu_{\lambda\lambda} =\begin{cases}\fr1{f_\lambda}\left(\fr1{f_\lambda} -\fr{f_{\lambda+1}}{(n_{\lambda+1}+1)f_K}\right),&{\rm if\ }1\leq a\leq \lambda;\\
\fr1{n_{\lambda+1}+1}\left(\fr1{n_{\lambda+1}+1}-\fr1{f_\lambda f_K}\right),&{\rm if\ }a=\lambda+1;\\
-\fr{f_{\lambda+1}}{(n_{\lambda+1}+1)f_\lambda f_K},&{\rm if\ }a\geq\lambda+2.
\end{cases}
\ee
Therefore, $\Gamma_{\lambda\lambda}^1,\,\cdots,\, \Gamma^{K-1}_{\lambda\lambda}$ satisfy the system same as in \eqref{gmij1} with the right hand side being replaced by
\begin{align}
&\left(\fr1{f_\lambda}\left(\fr1{f_\lambda} -\fr{f_{\lambda+1}}{(n_{\lambda+1}+1)f_K}\right),\cdots, \fr1{f_\lambda}\left(\fr1{f_\lambda} -\fr{f_{\lambda+1}}{(n_{\lambda+1}+1)f_K}\right),\right.\nnm\\
&\hs{3cm}\left. \fr1{n_{\lambda+1}+1}\left(\fr1{n_{\lambda+1}+1}-\fr1{f_\lambda f_K}\right),-\fr{f_{\lambda+1}}{(n_{\lambda+1}+1)f_\lambda f_K},\cdots,-\fr{f_{\lambda+1}}{(n_{\lambda+1}+1)f_\lambda f_K}\right)^t,\nnm
\end{align}
which is equivalent to the system \eqref{gmij2} with the right hand side being replaced by
$$
\left(0,\cdots,0, \fr1{f_\lambda^2}-\fr1{(n_{\lambda+1}+1)^2},\fr1{(n_{\lambda+1}+1)^2},0,\cdots, -\fr{f_{\lambda+1}}{(n_{\lambda+1}+1)f_\lambda f_K}\right)^t.
$$
Then it is easy to get $\Gamma_{\lambda\lambda}^1,\,\cdots,\, \Gamma^{K-1}_{\lambda\lambda}$ as in \eqref{gmlmdmu1}.
\endproof

{\lem\label{gamlmdia} It holds that
\begin{align}
&\Gamma^\mu_{\lambda\bar i_a}=0,\quad\Gamma^{\bar j_b}_{\lambda \bar i_a}=\delta^b_a\delta^{j_a}_{i_a}
\cdot\begin{cases} 0,&{\rm if\ }1\leq \lambda\leq a-2;\\
-\fr1{n_a+1},&{\rm if\ }\lambda=a-1;\\
\fr1{f_\lambda},&{\rm if\ }a\leq\lambda\leq K-1.
\end{cases}\label{gmlmdia}\\
&\Gamma^\lambda_{\bar i_a\bar j_b}=\Gamma^{\bar k_c}_{\bar i_a\bar j_b}=0,\quad{\rm if\ }a\neq b.\label{gmij3}
\end{align}}

\proof To prove \eqref{gmlmdia}, it suffices to consider
$$
\ppp{x}{t^\lambda}{u^{\bar i_a}}+L_1g_{\lambda \bar i_a} =\Gamma^\mu_{\lambda \bar i_a}\pp{x}{t^\mu}+\Gamma^{\bar j_b}_{\lambda \bar i_a}\pp{x}{u^{\bar j_b}},
$$
which is equivalent to that $\Gamma^{\bar j_b}_{\lambda \bar i_a}=0$ and
\begin{align}
&-\fr1{n_b+1}\Gamma^{b-1}_{\lambda\bar i_a} +\sum_{\mu=b}^{K-1}\fr1{f_\mu}\Gamma^\mu_{\lambda\bar i_a} =0,\label{gmlmdiamu}\\
&\Gamma^{\bar j_b}_{\lambda\bar i_a}=\delta^b_a\delta^{j_a}_{i_a}\begin{cases}0,&{\rm if\ }1\leq \lambda\leq a-2;\\
-\fr1{n_a+1},&{\rm if\ }\lambda=a-1;\\
\fr1{f_\lambda},&{\rm if\ }a\leq\lambda\leq K-1.
\end{cases}\label{gmlmdiajb}
\end{align}
This proves \eqref{gmlmdia}.

Moreover, if $a\neq b$, we have  that
$$
\ppp{x}{u^{\bar i_a}}{u^{\bar i_b}}+L_1g_{\bar i_a\bar j_b} =0,
$$
implying the vanishing of $\Gamma^\lambda_{\bar i_a\bar j_b}$ and $\Gamma^{\bar k_c}_{\bar i_a\bar j_b}$ for $a\neq b$.
\endproof

\newcommand{\oma}{\stx{a}{\omega}{}\!\!}
\newcommand{\omb}{\stx{b}{\omega}{}\!\!}
Let $\{\pp{x}{u^i},e_{n+1}\}$ and $\{\pp{x_a}{v^{i_a}},e^{(a)}_{n_a+1}\}$ ($a=1,\cdots,K$) be unimodular with $e_{n+1}$ and $e^{(a)}_{n_a+1}$ ($a=1,\cdots,K$) parallel to the corresponding affine normals. We denote by $\omega^j_i$, $\omega^{n+1}_{n+1}$ (resp. $\oma^{j_a}_{i_a}$, $\oma^{n_a+1}_{n_a+1}$) the induced affine connection form of the hyperbolic affine hypersphere $x$ (resp. $x_a$). Then we have (\cite{li-sim-zhao93}):
\be\label{fml9'}
\omega^{n+1}_{n+1}+\fr1{n+2}d\log H=0,\quad \oma^{n_a+1}_{n_a+1}+\fr1{n_a+2}d\log\Ha=0,\ a=1,\cdots,K.
\ee

Then by \eqref{fml9'} and Lemmas \ref{gamij}--\ref{gamlmdia} we have the following corollary:

{\cor\label{connform} The induced affine connection forms $\{\omega^j_i,\omega^{n+1}_{n+1}\}$ with respect to $\{\pp{x}{u^i},e_{n+1}\}$ are given by \begin{align}
\omega^\lambda_\lambda=&\sum_{\mu}\Gamma^\lambda_{\lambda\mu}dt^\mu +\sum_{a,i_a}\Gamma^\lambda_{\lambda\bar i_a}du^{\bar i_a}=\sum_{\lambda}\left(\fr1{f_\lambda}-\fr1{n_{\lambda+1}+1}\right)dt^\lambda +\sum_{\mu>\lambda}\fr1{f_\mu}dt^\mu;\label{omega1}\\
\omega^{\mu}_{\lambda}=&\sum_{\nu}\Gamma^{\mu}_{\lambda\nu}dt^\nu +\sum_{a,i_a}\Gamma^{\mu}_{\lambda\bar i_a}du^{\bar i_a}=
\fr{n_{\mu+1}}{f_{\mu+1}}\sum_{\nu}\fr{f_{\nu+1}}{(n_{\nu+1})f_\nu}dt^\nu, \quad{\rm for\ }\lambda<\mu;\label{omega2}\\
\omega^{\mu}_{\lambda}=&\sum_{\nu}\Gamma^{\mu}_{\lambda\nu}dt^\mu +\sum_{a,i_a}\Gamma^{\mu}_{\lambda\bar i_a}du^{\bar i_a}=\fr1{f_\lambda}dt^\mu,\quad{\rm for\ }\lambda>\mu;\label{omega3}\\
\omega^{\bar i_a}_{\lambda}=&\sum_{\mu}\Gamma^{\bar i_a}_{\lambda\mu}dt^\mu +\sum_{b,j_b}\Gamma^{\bar i_a}_{\lambda\bar j_b}du^{\bar j_b} =-\fr1{n_\lambda+1}\sum_{\bar j_{\lambda+1}}du^{\bar j_{\lambda+1}} +\fr1{f_\lambda}\sum_{a=1}^\lambda\sum_{\bar j_a} du^{\bar j_a};\label{omega4}\\
\omega^{\lambda}_{\bar i_a}=&\sum_{\mu}\Gamma^{\lambda}_{\bar i_a\mu}dt^\mu +\sum_{b,j_b}\Gamma^{\lambda}_{\bar i_a\bar j_b}du^{\bar j_b}=\fr{(n_{\lambda+1}+1)f_\lambda}{f_{\lambda+1}}\stx{\lambda+1}{L}\hs{-.3cm}_1 \hs{.1cm}\stx{\lambda+1}{g}\hs{-.36cm}_{i_{\lambda+1}j_{\lambda+1}}du^{\bar j_{\lambda+1}};\label{omega5}\\
\omega^{j_b}_{i_a}=&\sum_{\lambda}\Gamma^{\bar j_b}_{\bar i_a\lambda}dt^\lambda +\sum_{c,k_c}\Gamma^{\bar j_b}_{\bar i_a\bar k_c}du^{\bar k_c}=-\fr1{n_a+1}\delta^b_a\delta^{\bar j_b}_{\bar i_a}dt^{a-1} +\delta^b_a\sum_{\lambda\geq a}\fr1{f_\lambda}\delta^{j_b}_{i_a}dt^\lambda+\delta^b_a\oma^{j_a}_{i_a};\label{omega6}\\
\omega^{n+1}_{n+1}=&\sum_a\oma^{n_a+1}_{n_a+1}.\label{omega7}
\end{align}}

\subsection{The Fubini-Pick form}

Now we want to compute the Fubini-Pick form. To do this, the first step is to find the components $h_{ijk}$ which are defined by \eqref{hijk0}, namely
\be\label{hijk1}
\sum_{k}h_{ijk}\omega^k=dh_{ij}+h_{ij}\omega^{n+1}_{n+1}-\sum_{k}h_{kj}\omega^k_i -\sum_{k}h_{ik}\omega^k_j.
\ee
Similarly, we define $\ha_{i_aj_ak_a}$ for each $a$ by
\be\label{haijk}
\sum_{k_a}\ha_{i_aj_ak_a}\!\!\oma^{k_a}=d\ha_{i_aj_a}+\ha_{i_aj_a}\oma^{n_a+1}_{n_a+1} -\sum_{k_a}\ha_{k_aj_a}\oma^{k_a}_{i_a}-\sum_{k_a}\ha_{i_ak_a}\oma^{k_a}_{j_a}.
\ee

\newcommand{\tdca}{(n_a+1)\big(-\!\!\la\big)\prod_b\fr{c_b^{n_b+1}} {(n_b+1)\big(-\!\!\lb\big)}}
\newcommand{\cha}{\prod_a\fr{c_a^{n_a+1}\Ha^{\fr1{n_a+2}}} {(n_a+1)\big(-\!\!\la\big)}}
\newcommand{\chb}{\prod_b\fr{c_b^{n_b+1}\Hb^{\fr1{n_b+2}}} {(n_b+1)\big(-\!\!\lb\big)}}
\newcommand{\chc}{\prod_c\fr{c_c^{n_c+1}\Hc^{\fr1{n_c+2}}} {(n_c+1)\big(-\!\!\lc\big)}}
{\lem\label{hijk} The possibly non-zero components among $h_{ijk}$ are as follows:
\begin{align}
&h_{\lambda\lambda\lambda}=-2\left(\chb\right) \fr{f_{\lambda+1}}{(n_{\lambda+1}+1)f_\lambda}\left(\fr1{f_\lambda} -\fr1{n_{\lambda+1}+1}\right),\\
&h_{\lambda\lambda\mu}=-2\left(\chb\right) \fr{f_{\lambda+1}}{(n_{\lambda+1}+1)f_\lambda f_\mu},\quad\lambda<\mu,\\
&h_{\bar i_a\bar j_a\bar k_a}=\left(\chb\right)(n_a+1)\big(-\!\!\la\big)\Ha^{-\fr1{n_a+2}}\ha_{i_aj_ak_a},\\
&h_{\bar i_a\bar j_a\,a-1}=\fr{2}{(n_a+1)C}\left(\chb\right)g_{\bar i_a\bar j_a},\\
&h_{\bar i_a\bar j_a\lambda}=-\fr{2}{f_\lambda C}\left(\chb\right)g_{\bar i_a\bar j_a},\quad \lambda\geq a.
\end{align}}

\newcommand{\ca}{\prod_a\fr{c_a^{n_a+1}}{(n_a+1)\big(-\!\!\la)}}
\newcommand{\cb}{\prod_b\fr{c_b^{n_b+1}}{(n_b+1)\big(-\!\!\lb)}}
\newcommand{\cc}{\prod_c\fr{c_c^{n_c+1}}{(n_c+1)\big(-\!\!\lc)}}
\proof Firstly we compute $h_{\lambda\mu k}$. For this we use \eqref{hijk1}, \eqref{haijk} and Corollary \ref{connform} to find
\begin{align}
\sum_k h_{\lambda\mu k}\omega^k=&dh_{\lambda\mu}+h_{\lambda\mu}\omega^{n+1}_{n+1} -\sum_\nu h_{\nu\mu}\omega^\nu_{\lambda}-\sum_\nu h_{\lambda\nu}\omega^\nu_{\mu}\nnm\\
=&\ca\cdot\fr{f_{\lambda+1}}{(n_{\lambda+1}+1)f_\lambda} d\left(\prod_a\Ha^{\fr1{n_a+2}}\right)\delta_{\lambda\mu}\nnm\\
&\ +\cha\cdot\fr{f_{\lambda+1}}{(n_{\lambda+1}+1)f_\lambda}\delta_{\lambda\mu} \omega^{n+1}_{n+1}\nnm\\
&\ -\cha\cdot\fr{f_{\nu+1}}{(n_{\nu+1}+1)f_\nu}\delta_{\nu\mu} \omega^\nu_\lambda \nnm\\
&\ -\cha\cdot\fr{f_{\lambda+1}}{(n_{\lambda+1}+1) f_\lambda}\delta_{\lambda\nu} \omega^\nu_\mu\label{fml9}
\end{align}
But by \eqref{fml9'} and \eqref{H}
$$
\omega^{n+1}_{n+1}=-\fr1{n+2}d\log H=-\fr1{n+2}d\log\left(\prod_a\Ha^{\fr{n+2}{n_a+2}}\right) =-\prod_a\Ha^{-\fr1{n_a+2}}d\left(\prod_a \Ha^{\fr1{n_a+2}}\right).
$$
Thus
\be\label{fml9-1}
\left(\prod\Ha^{\fr1{n_a+2}}\right)\omega^{n+1}_{n+1}+d\prod_a \Ha^{\fr1{n_a+2}}=0.
\ee
This with \eqref{fml9} gives
\be\label{fml9-2}
\sum_k h_{\lambda\mu k}\omega^k=-\cha\left(\fr{f_{\mu+1}}{(n_{\mu+1}+1)f_\mu}\omega^\mu_\lambda
+\fr{f_{\lambda+1}}{(n_{\lambda+1}+1)f_\lambda}\omega^\lambda_\mu\right).
\ee

Case 1. $\lambda=\mu$. Then
\begin{align}
\sum_k h_{\lambda\lambda k}\omega^k=&-2\cha\fr{f_{\lambda+1}}{(n_{\lambda+1}+1) f_\lambda}\omega^\lambda_\lambda\nnm\\
=&-2\cha\fr{f_{\lambda+1}}{(n_{\lambda+1}+1)f_\lambda} \left(\left(\fr1{f_\lambda}-\fr1{n_{\lambda+1}+1}\right)dt^\lambda +\sum_{\mu=\lambda+1}^{K-1}\fr1{f_\mu}dt^\mu\right).\label{h-lmdlmdc}
\end{align}
It follows that $h_{\lambda\lambda \bar i_a}=0$ and
\bea
&&h_{\lambda\lambda\lambda}=-2\cha\cdot\fr{f_{\lambda+1}}{(n_{\lambda+1}+1)f_\lambda} \left(\fr1{f_\lambda}-\fr1{n_{\lambda+1}+1}\right),\\
&&h_{\lambda\lambda\mu}=-2\cha\cdot\fr{f_{\lambda+1}}{(n_{\lambda+1}+1)f_\lambda f_\mu}, \quad \lambda+1\leq\mu\leq K-1. \eea

Case 2. $\lambda<\mu$. In this case we have by \eqref{fml9-2} and Corollary \ref{connform}
\begin{align}
\sum_k h_{\lambda\mu k}\omega^k=&-\cha\cdot\left(\fr{f_{\mu+1}}{(n_{\mu+1}+1)f_\mu}\cdot \fr{(n_{\mu+1}+1)f_{\lambda+1}}{(n_{\lambda+1}+1)f_{\mu+1}f_\lambda}dt^\lambda +\fr{f_{\lambda+1}}{(n_{\lambda+1}+1)f_\lambda}\cdot \fr1{f_\mu}dt^\lambda\right)\nnm\\
=&-2\cha\cdot\left(\fr{f_{\lambda+1}}{(n_{\lambda+1}+1)f_\lambda f_{\mu}}dt^\lambda \right).
\end{align}
Thus $h_{\lambda\mu \bar i_a}=0$ for $\lambda<\mu$ and
\bea
&&h_{\lambda\mu\lambda}=-2\cha\cdot\fr{f_{\lambda+1}}{(n_{\lambda+1}+1)f_\lambda f_\mu},\\
&&h_{\lambda\mu\nu}=0,\quad{\rm for\ }\nu\neq\lambda.
\eea

Case 3. $\lambda>\mu$. In this case we have
\begin{align}
\sum_k h_{\lambda\mu k}\omega^k=&-\cha\cdot\left(\fr{f_{\mu+1}}{(n_{\mu+1}+1)f_\mu}\cdot \fr1{f_\lambda}dt^\mu+\fr{f_{\lambda+1}}{(n_{\lambda+1}+1)f_\lambda}\cdot \fr{(n_{\lambda+1}+1)f_{\mu+1}}{(n_{\mu+1}+1)f_{\lambda+1}f_\mu}dt^\mu \right)\nnm\\
=&-2\cha\cdot\left(\fr{f_{\mu+1}}{(n_{\mu+1}+1)f_\mu f_\lambda}dt^\mu\right).
\end{align}
Thus we have $h_{\lambda\mu \bar i_a}=0$ for $\lambda>\mu$ and
\bea
&&h_{\lambda\mu\mu}=-2\cha\cdot\fr{f_{\mu+1}}{(n_{\mu+1}+1)f_\mu f_\lambda},\\
&&h_{\lambda\mu\nu}=0,\quad{\rm for\ }\nu\neq\mu.
\eea

Secondly we compute $h_{\lambda\bar i_a k}$. This time we use the formula
\begin{align}
\sum_k h_{\lambda\bar i_a k}\omega^k=&dh_{\lambda\bar i_a}+h_{\lambda\bar i_a}\omega^{n+1}_{n+1} -\sum_k h_{k\bar i_a}\omega^k_{\lambda}-\sum_k h_{\lambda k}\omega^k_{\bar i_a}\nnm\\
=&-\sum_\mu h_{\mu\bar i_a}\omega^\mu_{\lambda} -\sum_{b,j_b}h_{\bar j_b\bar i_a}\omega^{\bar j_b}_{\lambda}-\sum_\mu h_{\lambda\mu}\omega^\mu_{\bar i_a} -\sum_{b,j_b}h_{\lambda\bar j_b}\omega^{\bar j_b}_{\bar i_a}=0.\label{fml10}
\end{align}
Consequently $h_{\lambda\bar i_a k}=0$, implying that $h_{\bar i_a\lambda k}=0$ by the symmetry.

Now we are left $h_{\bar i_a\bar j_b k}$ to compute. Note by \eqref{h ia jb} that
\begin{align}
h_{\bar i_a\bar j_b}=&\left(\cc\right)(n_a+1)\big(-\!\!\la\big) \left(\prod_c\Hc^{\fr1{n_c+2}}\right)\Ha^{-\fr1{n_a+2}} \ha_{i_aj_a}\delta_{ab}\nnm\\
=&\bar c_a\left(\prod_{c\neq a}\Hc^{\fr1{n_c+2}}\right)\ha_{i_aj_a}\delta_{ab},
\end{align}
where
$$
\bar c_a=\left(\cc\right)(n_a+1)\big(-\!\!\la\big).
$$

Case 1: $b=a$. We have
$$
dh_{\bar i_a\bar j_a}=\bar c_ad\left(\prod_{b\neq a}\Hb^{\fr1{n_b+2}}\right) \ha_{i_aj_a}+\bar c_a\left(\prod_{b\neq a}\Hb^{\fr1{n_b+2}}\right)d\ha_{i_aj_a}.
$$
It follows that
\begin{align}
&\sum_kh_{\bar i_a\bar j_a k}\omega^k\!=dh_{\bar i_a\bar j_a}+h_{\bar i_a\bar j_a}\omega^{n+1}_{n+1} -\sum_{k_a}h_{\bar k_a\bar j_a} \omega^{\bar k_a}_{\bar i_a} -\sum_{k_a}h_{\bar i_a \bar k_a}\omega^{\bar k_a}_{\bar j_a}\nnm\\
\!=&\bar c_ad\left(\prod_{b\neq a}\Hb^{\fr1{n_b+2}}\right)\ha_{i_aj_a} \!+\bar c_a\left(\prod_{b\neq a}\Hb^{\fr1{n_b+2}}\right)d\ha_{i_aj_a}\!+
\bar c_a\left(\prod_{b\neq a}\Hb^{\fr1{n_b+2}}\right) \ha_{i_aj_a}\sum_b\omb^{n_b+1}_{n_b+1}\nnm\\
&-\bar c_a\left(\prod_{b\neq a}\Hb^{\fr1{n_b+2}}\right)\left(\sum_{k_a}\ha_{k_aj_a}\omega^{\bar k_a}_{\bar i_b} +\sum_{k_a}\ha_{i_ak_a}\omega^{\bar k_a}_{\bar j_a}\right)\nnm\\
=&\bar c_ad\left(\prod_{b\neq a}\Hb^{\fr1{n_b+2}}\right)\ha_{i_aj_a} +\bar c_a\left(\prod_{b\neq a}\Hb^{\fr1{n_b+2}}\right) \left(d\ha_{i_aj_a}+\ha_{i_aj_a}\!\!\oma^{n_a+1}_{n_a+1}\right.\nnm\\ &\ \left.-\sum_{k_a}\ha_{k_aj_a}\!\!\oma^{k_a}_{i_a} -\sum_{k_a}\ha_{i_ak_a}\!\!\oma^{k_a}_{j_a}\right)+\bar c_a\left(\prod_{b\neq a}\Hb^{\fr1{n_b+2}}\right)\sum_{b\neq a}\omb^{n_b+1}_{n_b+1}\ha_{i_aj_a}\nnm\\
&\ -\bar c_a\left(\prod_{b\neq a}\Hb^{\fr1{n_b+2}}\right)\sum_{k_a}\left(-\fr1{n_a+1}\delta^{k_a}_{i_a}dt^{a-1} +\sum_{\lambda\geq a}\fr1{f_\lambda}\delta^{k_a}_{i_a}dt^\lambda\right)\ha_{k_aj_a}\nnm\\ &\ -\bar c_a\left(\prod_{b\neq a}\Hb^{\fr1{n_b+2}}\right)\sum_{k_a}\left(-\fr1{n_a+1}\delta^{k_a}_{i_a}dt^{a-1} +\sum_{\lambda\geq a}\fr1{f_\lambda}\delta^{k_a}_{i_a}dt^\lambda\right)\ha_{k_aj_a}\nnm\\
=&\bar c_a\left(d\left(\prod_{b\neq a}\Hb^{\fr1{n_b+2}}\right) +\left(\prod_{b\neq a}\Hb^{\fr1{n_b+2}}\right)\sum_{b\neq a}\omb^{n_b+1}_{n_b+1}\right)\ha_{i_aj_a} +\bar c_a\left(\prod_{b\neq a}\Hb^{\fr1{n_b+2}}\right)\sum_{k_a}\ha_{i_aj_ak_a}\oma^{k_a}\nnm\\
&-\bar c_a\left(\prod_{b\neq a}\Hb^{\fr1{n_b+2}}\right) \left(-\fr2{n_a+1}dt^{a-1}+\sum_{\lambda\geq a}\fr2{f_\lambda}dt^\lambda\right)\ha_{i_aj_a}\nnm\\
=&\bar c_a\left(\prod_b\Hb^{\fr1{n_b+2}}\right)\Ha^{-\fr1{n_a+2}} \sum_{k_a}\ha_{i_aj_ak_a}\oma^{k_a}+\fr2{n_a+1}\bar c_a\left(\prod_b\Hb^{\fr1{n_b+2}}\right) \Ha^{-\fr1{n_a+2}}\ha_{i_aj_a}dt^{a-1}\nnm\\
&-\bar c_a\left(\prod_b\Hb^{\fr1{n_b+2}}\right)\Ha^{-\fr1{n_a+2}}\ha_{i_aj_a}\sum_{\lambda\geq a}\fr2{f_\lambda}dt^\lambda.\label{fml11}
\end{align}
Therefore
\begin{align}
h_{\bar i_a\bar j_a\,a-1}=&\fr2{n_a+1}\left(\chb\right)(n_a+1)\big(-\!\!\la\big)\Ha^{-\fr1{n_a+2}} \ha_{i_aj_a}\nnm\\
=&\fr2{(n_a+1)C}\left(\chb\right)g_{\bar i_a\bar j_a},\nnm\\
h_{\bar i_a\bar j_a\,\lambda}=&-\fr2{f_\lambda C} \left(\chb\right)g_{\bar i_a\bar j_a},\quad a\leq\lambda\leq K-1,\nnm\\
h_{\bar i_a\bar j_a\bar k_b}=&\bar c_a\left(\prod_c\Hc^{\fr1{n_c+2}}\right)\Ha^{-\fr1{n_a+2}} \ha_{i_aj_ak_a}\delta_{ab}\nnm\\ =&\left(\chc\right)(n_a+1)\big(-\!\!\la\big)\Ha^{-\fr1{n_a+2}}\ha_{i_aj_ak_a}\delta_{ab}.
\end{align}

Case 2: $a\neq b$. By the fact that $h_{i_aj_b}=0$ (see \eqref{h ia jb}) and $\omega^{\bar j_b}_{\bar i_a}=0$ (see \eqref{omega6}), we easily find that
$h_{\bar i_a\bar j_b k}=0$ for all $k=1,\cdots,n$.
\endproof

\newcommand{\Aa}{\stx{a}{A}{}\!\!}
\newcommand{\Aalp}{\stx{\alpha}{A}{}\!\!}
Now let $\Aa_{i_aj_ak_a}$ be the components of the Fubini-Pick form of the hyperbolic affine hypersphere, $a=1,\cdots,K$. Then we have

{\prop\label{fubpicform} The possibly non-zero components of the Fubini-Pick form of the hyperbolic affine hypersphere $x$ are as follows:
\bea
&A_{\lambda\lambda\lambda}=\fr{f_{\lambda+1}C}{(n_{\lambda+1}+1) f_\lambda}\left(\fr1{f_\lambda} -\fr1{n_{\lambda+1}+1}\right),\\
&A_{\lambda\lambda\mu}=\fr{f_{\lambda+1}C}{(n_{\lambda+1}+1)
f_\lambda f_\mu},\quad\lambda<\mu,\\
&A_{\bar i_a\bar j_a\,a-1}=-\fr1{n_a+1}g_{\bar i_a\bar j_a}=-\!\la C\ga_{i_aj_a},\\
&A_{\bar i_a\bar j_a\lambda}=\fr1{f_\lambda}g_{\bar i_a\bar j_a}=\fr{(n_a+1)\big(-\!\!\la\big)C}{f_\lambda}\ \ga_{i_aj_a},\quad \lambda\geq a,\\
&A_{\bar i_a\bar j_a\bar k_a}=(n_a+1)\big(-\!\!\la\big)C\Aa_{i_aj_ak_a}.
\eea}

\proof We shall use formula \eqref{hijktoaijk} to find the Fubini-Pick form of $x$. Note that $f_K+1=n+2$. Thus, by \eqref{H} and the definition of the constant $C$ in \eqref{l1c}, it holds that
\be
H^{-\fr1{n+2}}\cdot\cha=\h\cdot\cha=C.
\ee
Then the present proposition follows readily from Lemma \ref{hijk}.
\endproof

The Levi-Civita connection $\ol\omega^j_i=\ol\Gamma^j_{ik}du^k$ of the Berwald-Blaschke $g$ of $x$ is defined by
$$
\ol\Gamma^j_{ik}=\fr12g^{jl}\left(\pp{g_{il}}{u^k}+\pp{g_{kl}}{u^i} -\pp{g_{ik}}{u^j}\right).
$$
\newcommand{\olomea}{\stx{a}{\ol\omega}{}\!\!}
\newcommand{\olomeb}{\stx{b}{\ol\omega}{}\!\!}
\newcommand{\olgma}{\stx{a}{\ol\Gamma}{}\!\!\!}
\newcommand{\olgmb}{\stx{b}{\ol\Gamma}{}\!\!\!}
Similarly defined, for each $a$, is the Levi-Civita connection
$\olomea^{j_a}_{i_a}=\olgma^{j_a}_{i_ak_a}dv^{k_a}$ of the metric $\ga$ of $x_a$. Those formulas together with \eqref{g-lmdmu}--\eqref{g-lmdia} show that
\be\label{levi}
\ol\omega^\mu_\lambda=\ol\omega^{\bar i_a}_\lambda=\ol\omega^\lambda_{\bar i_a}=0,\quad \ol\omega^{\bar j_b}_{\bar i_a}=\delta^b_a\olomea^{j_a}_{i_a}.
\ee

An application of Proposition \ref{fubpicform} and \eqref{levi} proves the following

{\prop\label{paral pic} ({\rm cf.} \cite{dil-vra94}, Proposition 2) With respect to the Levi-Civita connections, the Calabi composition $x$ is of parallel Fubini-Pick form if and only if each factor $x_a$ is.}

\section{Extending to the zero dimensional factors---a more general formula}

Now we extend the Calabi composition discussed above to
some more general setting.

For each point in $\bbr$ different from the origin or, equivalently, for each nonzero number $c$, we can view it as a zero-dimensional ``hyperbolic affine hypersphere'' immersed in $\bbr$ with the affine mean curvature being $-1$. In this sense, we may allow some or all of the factors $x_a$ of the Calabi composition to be replaced by some given points or some given positive numbers. In fact, a step-by-step examination shows that all the argument and the computations in the previous subsections apply in this general situation. Without loss of generality, we can normalize these given $0$-dimensional factors to be $1$. Thus all the conclusions we have obtained previously remain true if we allow some or all of $n_a$'s to be zero. For example, Theorem \ref{multicomp} can be refreshed as follows:

{\thm\label{general sense} Let $r,s$ be two nonnegative integers with $K:=r+s\geq 2$ and $x_\alpha:M^{n_\alpha}\to\bbr^{n_\alpha+1}$, $1\leq \alpha\leq s$, be hyperbolic affine hyperspheres of dimension $n_\alpha>0$ with affine mean curvatures $\stx{\alpha}{L}\!\!_1$ and with the origin their common affine center. Then, for any $K$ positive numbers $c_1,\cdots,c_K$, we have a new hyperbolic affine hypersphere $x:M^n\to\bbr^{n+1}$ with the affine mean curvature
\be\label{newl1c}
L_1=-\fr1{(n+1)C},\quad C:=\left(\fr1{n+1}\prod_{a=1}^r c_a^2\cdot\prod_{\alpha=1}^s\fr{c_{r+\alpha}^{2(n_\alpha+1)}} {(n_\alpha+1)^{n_\alpha+1}(-\!\!\stx{\alpha}{L}_1)^{n_\alpha+2}}\right)^{\fr1{n+2}},
\ee
where $n=\sum_\alpha n_\alpha+K-1$, $M^n=\bbr^{K-1}\times M^{n_1}\times\cdots\times M^{n_s}$ and
\begin{align}
x(t^1,&\cdots,t^{K-1},p_1,\cdots,p_s):=(c_1e_1,\cdots, c_re_r,c_{r+1}e_{r+1}x_1(p_1),\cdots,c_Ke_Kx_s(p_s)),\nnm\\&\hs{1cm}\forall (t^1,\cdots,t^{K-1},p_1,\cdots,p_s)\in M^n.\label{mulpro2}
\end{align}

Moreover, for given positive numbers $c_1,\cdots,c_K$, there exits some $c>0$ and $c'>0$ such that
the following three hyperbolic affine hyperspheres
\bea &x:=(c_1e_1,\cdots, c_re_r,c_{r+1}e_{r+1}x_1,\cdots,c_Ke_sx_s),\nnm\\
&\bar x:=c(e_1,\cdots, e_r,e_{r+1}x_1,\cdots,e_sx_s),\nnm\\
&\td x:=(e_1,\cdots, e_r,e_{r+1}x_1,\cdots,ce_sx_s)\nnm
\eea
are equiaffine equivalent to each other.}

{\dfn\label{df2}\rm The hyperbolic affine hypersphere $x$ is called the Calabi composition of $r$ points and $s$ hyperbolic affine hyperspheres.}

\rmk\rm Two special cases of the above proposition when $r=0,s=2$ and $r=s=1$, respectively, are discussed in \cite{dil-vra94} and \cite{hu-li-vra08}.

Now let $x:M^n\to \bbr^{n+1}$ be a Calabi composition of $r$ points and $s$ ($K=r+s\geq 2$) hyperbolic affine hyperspheres $x_\alpha:M^{n_\alpha}\to\bbr^{n_\alpha+1}$, $1\leq \alpha\leq s$, and $g$ the affine metric of $x$. For convenience we make the following convention:
$$
1\leq\alpha,\beta,\gamma\leq s,\quad \td\alpha=\alpha+r,\ \td\beta=\beta+r,\ \td\gamma=\gamma+r.
$$

Let $\{v^{i_\alpha}_\alpha; \ i_\alpha=1,\cdots,n_\alpha\}$ be the local coordinate system of $M_\alpha$, $\alpha=1,\cdots,s$. Note that
$$
f_a=\begin{cases} a,&1\leq a\leq r;\\ \sum_{\beta\leq \alpha}n_\beta+\td{\alpha},&r+1\leq a=\td\alpha\leq r+s.
\end{cases}
$$
Let $C$ be the constant given by \eqref{newl1c}. Then from \eqref{g-lmdmu}, \eqref{g-iajb}, \eqref{g-lmdia} and Proposition \ref{fubpicform} follows easily the following

{\prop\label{corr0} The affine metric $g$, the affine mean curvature $L_1$ and the possibly nonzero components of the Fubini-Pick form $A$ of the Calabi composition $x$ of $r$ points and $s$ hyperbolic affine hyperspheres $x_\alpha:M_\alpha\to\bbr^{n_\alpha+1}$, $\alpha=1,\cdots,s$, are given as follows:
\begin{align}
&g_{\lambda\mu}=\begin{cases}\displaystyle \fr{\lambda+1}{\lambda}C\delta_{\lambda\mu},&1\leq\lambda\leq r-1;\\
\displaystyle\fr{n_1+r+1}{r(n_1+1)}C\delta_{r\mu},&\lambda=r;\\
\displaystyle\fr{\sum_{\beta\leq\alpha+1}n_\beta+\td{\alpha}+1}  {(n_\alpha+1)(\sum_{\beta\leq \alpha}n_\beta+\td{\alpha})}C\delta_{\lambda\mu}, &r+1\leq\lambda=\td\alpha\leq r+s-1.
\end{cases}
\\
&g_{\bar i_{\td\alpha}\bar j_{\td\beta}}=(n_\alpha+1) (-\!\!\lalp)C\galp_{i_{\alpha}j_{\alpha}}\delta_{\alpha\beta},
\quad
g_{\lambda\bar i_{\td\alpha}}=0.\label{g-gen}
\\
&A_{\lambda\lambda\lambda} =\begin{cases}\displaystyle \fr{1-\lambda^2}{\lambda^2}C,&1\leq\lambda\leq r-1,\\
\displaystyle\left(\fr1{r^2}-\fr1{(n_1+1)^2}\right)C,&\lambda=r,\\
\displaystyle\fr{(\sum_{\beta\leq\alpha+1}n_\beta+\td{\alpha}+1)C} {(n_{\alpha+1}+1)(\sum_{\beta\leq \alpha}n_\beta+\td{\alpha})}\left(\fr1{\sum_{\beta\leq \alpha}n_\beta+\td{\alpha}}-\fr1{n_{\alpha+1}+1}\right), &r+1\leq\lambda=\td\alpha\leq r+s-1.
\end{cases}
\\
&A_{\lambda\lambda\mu} =\begin{cases}\displaystyle \fr{\lambda+1}{\lambda\mu}C,&1\leq\lambda<\mu\leq r,\\
\displaystyle\fr{(\lambda+1)C}{\lambda(\sum_{\beta\leq \alpha}n_\beta+\td{\alpha})},&1\leq\lambda\leq r-1, \mu=\td\alpha,\\
\displaystyle\fr{(n_1+r+1)C}{r(\sum_{\beta\leq \alpha}n_\beta+\td{\alpha})},&\lambda=r,\ \mu=\td\alpha,\\
\displaystyle\fr{(\sum_{\gamma\leq \alpha+1}n_\gamma+\td{\alpha}+1)C} {(n_{\alpha+1}+1)(\sum_{\gamma\leq \alpha}n_\gamma+\td{\alpha})(\sum_{\gamma\leq \beta}n_\gamma+\td{\beta})},&r+1\leq\lambda=\td\alpha<\mu=\td\beta\leq r+s-1.
\end{cases}
\\
&A_{\bar i_{\td\alpha}\bar j_{\td\alpha}\,{\td\alpha}-1} =-\fr1{n_\alpha+1}g_{\bar i_{\td\alpha}\bar j_{\td\alpha}}=-\!\lalp C\galp_{i_\alpha j_\alpha},\\
&A_{\bar i_{\td\alpha}\bar j_{\td\alpha}\td\beta}=\fr1{\sum_{\gamma\leq \beta}n_\gamma+\td{\beta}}g_{\bar i_{\td\alpha}\bar j_{\td\alpha}}=\fr{(n_\alpha+1)\big(-\!\!\lalp\big)C}{\sum_{\gamma\leq \beta}n_\gamma+\td{\beta}}\ \galp_{i_\alpha j_\alpha},\quad \beta\geq \alpha,\\
&A_{\bar i_{\td\alpha}\bar j_{\td\alpha}\bar k_{\td\alpha}}=(n_\alpha+1)\big(-\!\!\lalp\big)C\Aalp_{i_\alpha j_\alpha k_\alpha},
\end{align}
where $\stx{\alpha}{L}_1$, $\stx{\alpha}{g}$ and $\stx{\alpha}A$ are the affine mean curvature, the affine metric and the Fubini-Pick form of $x_\alpha$, $\alpha=1,\cdots,s$.}

By restrictions, $g$ defines a flat metric $g_0$ on $\bbr^{K-1}$ with matrix $(g_{\lambda\mu})$ and, for each $\alpha$, a metric $g_\alpha$ on $M_\alpha$ with matrix $\big(g^\alpha_{i_\alpha j_\alpha}\big)=\big(g_{\bar i_{\td\alpha}\bar j_{\td\alpha}}\big)$ and inverse matrix $\big(g^{i_\alpha j_\alpha}_\alpha\big)$, which is conformal to the original metric $\stx{\alpha}{g}$, or more precisely, \be\label{g_alpha}g_\alpha=(n_\alpha+1) \big(-\!\!\stx{\alpha}{L}_1)C\stx{\alpha}{g}.\ee

\expl\label{expl}\rm Given a positive number $C_0$, let $x_0:\bbr^{n_0}\to \bbr^{n_0+1}$ be the well known flat hyperbolic affine hypersphere of dimension $n_0$ which is defined by
$$
x^1\cdots x^{n_0} x^{n_0+1}=C_0,\quad x^1>0,\cdots,x^{n_0+1}>0.
$$
Then it is not hard to see that $x_0$ is the Calabi composition of $n_0+1$ points. In fact, we can write for example
$$
x_0=(e_1,\cdots,e_{n_0},C_0e_{n_0+1}).
$$
Then by Corollary \ref{corr0} the affine metric $g_0$, the affine mean curvature $\stx{0}{L}_1$ and the Fubini-Pick form $\stx{0}{A}$ of $x_0$ are respectively given by (cf. \cite{li-sim-zhao93})
\begin{align}
\stx{0}{g}_{\lambda\mu}=& \fr{\lambda+1}{\lambda}\left(\fr{C_0^2}{n_0+1}\right)^{\fr1{n_0+2}} \delta_{\lambda\mu},\label{gofexpl}\\
\stx{0}{L}_1=&-\fr1{(n_0+1)C} =-(n_0+1)^{-\fr{n_0+1}{n_0+2}}C_0^{-\fr2{n_0+2}},\label{l1ofexpl}\\
\stx{0}{A}_{\lambda\mu\nu}=&\begin{cases}-\fr{\lambda^2-1}{\lambda^2} \left(\fr{C_0^2}{n_0+1}\right)^{\fr1{n_0+2}}, &{\rm if\ }\lambda=\mu=\nu;\\
\fr{\lambda+1}{\lambda\nu}\left(\fr{C_0^2}{n_0+1}\right)^{\fr1{n_0+2}}, &{\rm if\ }\lambda=\mu<\nu;\\
0, &{\rm otherwise.}
\end{cases}\label{aofexpl}
\end{align}
Thus the Pick invariant of $x_0$ is
\be\label{jofexpl}
\stx{0}{J}=\fr1{n_0(n_0-1)}\stx{0}{g}{}\!\!^{\lambda_1\lambda_2} \stx{0}{g}{}\!\!^{\mu_1\mu_2} \stx{0}{g}{}\!\!^{\nu_1\nu_2} \stx{0}{A}_{\lambda_1\mu_1\nu_1}\stx{0}{A}_{\lambda_2\mu_2\nu_2} =(n_0+1)^{-\fr{n_0+1}{n_0+2}}C_0^{-\fr2{n_0+2}}=-\!\!\stx{0}{L}_1.
\ee

\rmk\rm Since, by Propositions \ref{comml}, \ref{assl}, Calabi composition of hyperbolic affine hyperspheres is essentially independent of the order of its factors, we conclude that the Calabi composition of $r$$(\geq 2)$ points and $s$ hyperbolic affine hyperspheres is in fact the composition of the flat hyperbolic affine hypersphere $x_0$ in Example \ref{expl} of dimension $n_0:=r-1$ with other $s$ hyperbolic affine hyperspheres.

As an application of Proposition \ref{corr0}, we prove the following result:

{\thm\label{added} Let $x:M^n\to \bbr^{n+1}$ be a Calabi composition of $r$ points and $s$ hyperbolic affine hyperspheres $x_\alpha:M^{n_\alpha}\to\bbr^{n_\alpha+1}$, $1\leq \alpha\leq s$. Then $x$ is symmetric if and only if
each positive dimensional factor $x_\alpha$ is symmetric.}

\begin{proof} By using Proposition \ref{corr0}, it suffices to prove the following

{\prop\label{sym para1}
A nondegenerate hypersurface $x:M^n\to\bbr^{n+1}$ is locally affine symmetric if and only if $x$ is of parallel Fubini-Pick form $A$.}

{\it Proof of Propostion \ref{sym para1}}

First we suppose that the Fubini-Pick form $A$ of $x$ is parallel. Then by \cite{bok-nom-sim90}, $x$ must be an affine hypersphere. It then follows from \eqref{gaus_af sph} that the affine metric $g$ must be locally symmetric. Thus locally we can write $M^n=G/K$ and the canonical decomposition of the corresponding orthogonal symmetric pair $(\mathfrak{g},\mathfrak{k})$ is written as
${\mathfrak g}=\mathfrak{k}+\mathfrak{m}$ where the vector space $\mathfrak m$ is identified with $T_oM$. Here $o\in M^n$ is the base point given by $o=eK$ with $e$ the identity of $G$.
Note that, for all $X,Y_i\in {\mathfrak m}=T_oM$, $i=1,2,3$, the vector field $Y_i(t):=L_{\exp (tX)*}(Y_i)$ is the parallel translation of $Y_i$ along the geodesic $\gamma(t):=_{\exp (tX)}\!\!K$ (see, for example, \cite{hel01}).
Consequently we have
\begin{align}
&\dd{}{t}((L_{\exp (tX)}^*A)(Y_1,Y_2,Y_3))\nnm\\
=&\dd{}{t}(A_{\exp (tX)K}(L_{\exp (tX)*}(Y_1),
L_{\exp (tX)*}(Y_2),L_{\exp (tX)*}(Y_3)))\nnm\\
=&(\hat\nabla_{\gamma'(t)}A)(Y_1(t),Y_2(t),Y_3(t))=0,
\end{align}
where $\hat{\nabla}$ is the Levi-Civita connection of the metric $g$.
It follows that
\be\label{2.19-0}
A_{\exp (tX)K}(L_{\exp (tX)*}(Y_1),
L_{\exp (tX)*}(Y_2),L_{\exp (tX)*}(Y_3))
\ee
is constant with respect to the parameter $t$ and thus $A$ is $G$-invariant.

Conversely, we suppose that $M^n=G/K$ locally for some symmetric pair $(G,K)$ and that $A$ is $G$-invariant. Then for any $X,Y_i\in {\mathfrak m}=T_oM$, $i=1,2,3$, the function
\eqref{2.19-0}
is again a constant along the geodesic $\gamma(t)$.

Therefore,
$$
(\hat\nabla_XA)(Y_1,Y_2,Y_3)=\left.\dd{}{t}\right|_{t=0} A_{\gamma(t)}(Y_1(t),Y_2(t),Y_3(t))=0,
$$
where we have once again used the fact that each $Y_i(t)$ is parallel along the geodesic $\gamma(t)$.
\end{proof}

To end this section we list some properties of the Calabi composition of points and hyperbolic affine hyperspheres.

Write $M_0=\bbr^{K-1}$. Then, with respect to the affine metric $g$ on $M^n$, the Fubini-Pick form $A$ can be identified with a $TM^n$-valued symmetric $2$-form $ A:TM^n\times TM^n\to TM^n$. For each ordered triple $\alpha,\beta,\gamma\in\{0,1,\cdots,s\}$, $ A$ defines one $TM_\gamma$-valued bilinear map $ A^{\gamma}_{\alpha\beta}:TM_{\alpha}\times TM_{\beta}\to TM_\gamma$, which is the $TM_\gamma$-component of $ A_{\alpha\beta}$, the restriction of $ A$ to $TM_{\alpha}\times TM_{\beta}$. For $\alpha=1,\cdots,s$, define
\be\label{Halpha}
H_\alpha=\fr1{n_{\alpha}}\tr_{g_\alpha} A^0_{\alpha\alpha}\equiv \fr1{n_{\alpha}}g^{i_\alpha j_\alpha}_\alpha A^0_{\alpha\alpha}\left(\pp{}{v^{i_{\alpha}}_{\alpha}}, \pp{}{v^{j_{\alpha}}_{\alpha}}\right),
\ee
where the metric $g_\alpha$ is given by \eqref{g_alpha}.

\begin{prop}\label{properties} Let $x:M^n\to \bbr^{n+1}$ be the Calabi composition of $r$ points and $s$ hyperbolic affine hyperspheres and $g$ the affine metric of $x$. Then

$(1)$ The Riemannian manifold $M^n\equiv (M^n,g)$ is reducible, that is
\be\label{derham}
(M^n,g)=\bbr^q\times(M_1,g_1)\times\cdots\times(M_s,g_s),\quad q+s\geq 2;\ee

$(2)$ There must be a positive dimensional Euclidean factor $\bbr^q$ in the de Rham decomposition \eqref{derham} of $M^n$, that is, $q>0$;

$(3)$ $q\geq s-1$ with the equality holding if and only if $r=0$;

$(4)$ $ A^\gamma_{\alpha\beta}\equiv 0$ if
$(\alpha,\beta,\gamma)$ is not one of the following triples: $(0,0,0)$, $(\alpha,\alpha,0)$, $(\alpha,0,\alpha)$, $(0,\alpha,\alpha)$ or $(\alpha,\alpha,\alpha)$.

$(5)$ For any $p=(p_0,p_1,\cdots,p_s)\in M^n$ and each $\alpha=1,\cdots,s$, it holds that
\begin{align}\label{added f0}
A^0_{\alpha\alpha}(R^{M_\alpha}(X_\alpha,Y_\alpha)&Z_\alpha,W_\alpha)+A^0_{\alpha\alpha}(Z_\alpha,R^{M_\alpha}(X_\alpha,Y_\alpha)W_\alpha)=0.\nnm\\
&\forall X_\alpha,Y_\alpha,Z_\alpha,W_\alpha\in T_{p_\alpha}M_\alpha
\end{align}
which is equivalent to that the holonomy algebra ${\mathfrak h}_\alpha$ acts on $A^0_{\alpha\alpha}$ trivially, that is
${\mathfrak h}_\alpha\cdot A^0_{\alpha\alpha}=0$.

$(6)$ The vector-valued functions $H_\alpha$, $\alpha=1,\cdots,s$, defined by \eqref{Halpha} satisfy the following qualities:
\begin{align}
&H_\alpha=-\fr{f_{\td\alpha-1}}{f_{\td\alpha}C}\pp{}{t^{\td\alpha-1}} +\fr1C\sum_{s-1\geq\beta\geq\alpha}\fr{n_{\beta+1}+1} {f_{\td\beta+1}}\pp{}{t^{\td\beta}},\label{Halpha1}\\
&g(H_\alpha,H_\alpha)=C^{-1}\left(\fr1{n_{\alpha}+1}-\fr1{f_K}\right) =\fr{n-n_\alpha}{n_\alpha+1}(-L_1),\\
&g(H_\alpha,H_\beta)=L_1\quad{\rm for\ }\alpha\neq\beta;
\end{align}

$(7)$ $A^{\alpha}_{\alpha\alpha}$ is identical to the $TM_\alpha$-valued symmetric bilinear form defined by the Fubini-Pick form $\stx{\alpha}{A}$ of $x_\alpha$.
\end{prop}

In fact, (1), (2) are trivial, and (3) holds since $q=K-1=r+s-1$; (4) follows directly from Proposition \ref{corr0}; (5), (6) and (7) can be verified via a direct computation, where the equality \eqref{newl1c} and the fact that $f_K=n+1$ are needed.

\rmk\rm Another application of Proposition \ref{corr0} is to give a new characterization of the Calabi compostion. We believe that a locally strongly convex affine hypersurface $x:M^n\to \bbr^{n+1}$ is locally the Calabi composition of some points and hyperbolic affine hyperspheres if and only if the above conditions (1), (4) and (5) hold. This new characterization will appear in a forthcoming paper.

\end{document}